\documentclass[11pt,reqno]{amsart}
\usepackage[left=1in,top=1in,right=1in,bottom=1in]{geometry}
\usepackage{amsfonts}
\usepackage{enumerate,subfigure}
\usepackage{amssymb,latexsym, amsmath,graphicx}
\usepackage{epstopdf}
\usepackage{framed}
\usepackage{hyperref}
\usepackage{cite}

\usepackage[dvipsnames]{xcolor}

\numberwithin{equation}{section}

\usepackage{multirow}

\DeclareMathOperator{\R}{{\mathbb R}}
\DeclareMathOperator{\C}{{\mathbb C}}
\DeclareMathOperator{\by}{\times}
\DeclareMathOperator{\bndry}{\partial\Omega}
\DeclareMathOperator{\texp}{\mathbf{t}^{\mbox{\tiny{\textbf{exp}}}}}
\DeclareMathOperator{\qexp}{\mathbf{q}^{\mbox{\tiny{\textbf{exp}}}}}
\DeclareMathOperator{\sigcal}{\gamma^{\mbox{\tiny{\textbf{CAL}}}}}
\DeclareMathOperator{\sigexpCal}{\gamma^{\mbox{\tiny{\textbf{exp}}}}_{\mbox{\tiny{\textbf{CAL}}}}}
\DeclareMathOperator{\texpCAL}{\mathbf{t}^{\mbox{\tiny{\textbf{exp}}}}_{\mbox{\tiny{\textbf{CAL}}}}}
\DeclareMathOperator{\sigexp}{\gamma^{\mbox{\tiny{\textbf{exp}}}}}

\DeclareMathOperator{\Fhat}{\hat{F}}
\DeclareMathOperator{\gBest}{\gamma_{\mbox{\tiny{best}}}}
\DeclareMathOperator{\gammaB}{\gamma_b}

\newtheorem{prop}{Proposition}

\makeatletter
\def\thickhline{%
  \noalign{\ifnum0=`}\fi\hrule \@height \thickarrayrulewidth \futurelet
   \reserved@a\@xthickhline}
\def\@xthickhline{\ifx\reserved@a\thickhline
               \vskip\doublerulesep
               \vskip-\thickarrayrulewidth
             \fi
      \ifnum0=`{\fi}}
\makeatother

\newlength{\thickarrayrulewidth}
\usepackage{boldline}

\definecolor{brown}{rgb}{0.7,0.3,0}

\begin{document}

\title[3D CGO-Based EIT with Electrode Data]{3D EIT Reconstructions from Electrode Data using Direct Inversion D-bar and Calder\'on Methods}

\author{S.~J. Hamilton, D. Isaacson, ~V. Kolehmainen, ~P.~A. Muller, ~J. Toivanen, and ~P.~F. Bray}


\thanks{S.~J. Hamilton is with the Department of Mathematical and Statistical Sciences; Marquette University, Milwaukee, WI 53233 USA,  email: \texttt{sarah.hamilton@marquette.edu}}
\thanks{D. Isaacson is with the Department of Mathematical Sciences; Rensselaer Polytechnic Institute, Troy, NY 12180 USA,  email: \texttt{isaacd@rpi.edu}}
\thanks{V. Kolehmainen is with the Department of Applied Physics, University of Eastern Finland, FI-70210 Kuopio, Finland, email: \texttt{ville.kolehmainen@uef.fi}}
\thanks{P.~A. Muller* is with the Department of Mathematics \& Statistics; Villanova University, Villanova, PA 19085 USA,  email: \texttt{peter.muller@villanova.edu}, (corresponding author)}
\thanks{J. Toivanen is with the Department of Applied Physics, University of Eastern Finland, FI-70210 Kuopio, Finland, email: \texttt{jussi.toivanen@uef.fi}}
\thanks{P.~F. Bray is with the Department of Mathematics; Drexel University, Philadelphia, PA 19104 USA, email: \texttt{pfb25@drexel.edu}}

\begin{abstract}
The first numerical implementation of a D-bar method in 3D using electrode data is presented.  Results are compared to Calder\'on's method as well as more common TV and smoothness regularization-based methods.  D-bar methods are based on tailor-made non-linear Fourier transforms involving the measured current and voltage data.  Low-pass filtering in the non-linear Fourier domain is used to stabilize the reconstruction process.  D-bar methods have shown great promise in 2D for providing robust real-time absolute and time-difference conductivity reconstructions but have yet to be used on practical electrode data in 3D, until now.  Results are presented for simulated data for conductivity and permittivity with disjoint non-radially symmetric targets on spherical domains and noisy voltage data.  The 3D D-bar and Calder\'on methods are demonstrated to provide comparable quality to their 2D CGO counterparts, and hold promise for real-time reconstructions.
\end{abstract}
\maketitle 


\section{Introduction}\label{sec:Intro}

\begin{figure}[h!]
\centering
\includegraphics[width=280pt]{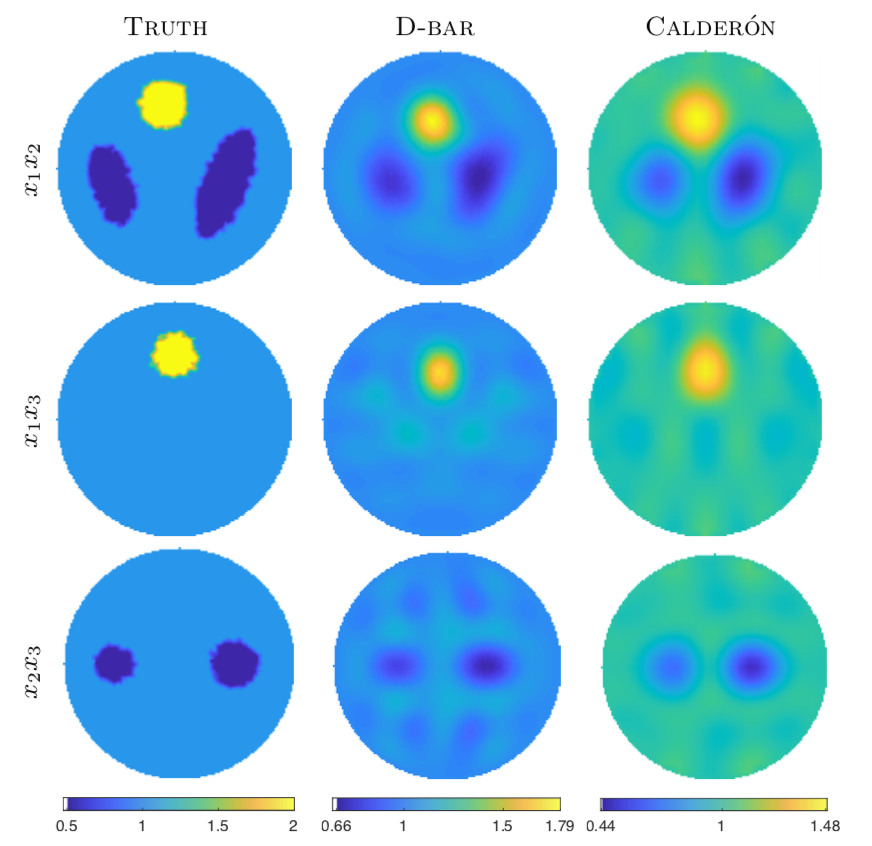}
\caption{\label{fig:Hammer} Demonstration of the 3D D-bar Method and Calder\'on's Method on the `Heart and Lungs' phantom T2-B using simulated electrode data.  The 2D cross-sectional slices above show that the conductive heart is correctly visible in the $x_1x_3$ plane but absent from the $x_2x_3$ plane.  Similarly for the lungs in the $x_2x_3$ plane vs the $x_1x_3$ plane.}
\end{figure}

In this manuscript, the first reconstructions using electrode data, instead of continuum boundary data, are presented for the D-bar method in three-dimensions.  The first comparison of D-bar, Calder\'on, and traditional regularized non-linear least squares based methods is shown.  Results in 3D are similar to those from 2D CGO-based methods in both resolution and recovered conductivity values. Figure~\ref{fig:Hammer} shows reconstructions for the `heart and lungs' phantom in 3D.

In Electrical Impedance Tomography (EIT), harmless currents are applied on electrodes affixed to the surface of a body and the resulting voltages are measured.  Using these pairs of surface electrical measurements, the task is to recover the point-wise conductivity and permittivity distributions {\it inside} the domain to `see inside'.  This is done by solving a severely ill-posed, non-linear inverse problem derived from Maxwell's equations and exploits the fact that internal tissues and inclusions have distinctive conductivity properties.  Medical applications of EIT include, but are not limited to, pulmonary function monitoring, detection of cancerous tissue, and imaging of stroke see e.g. \cite{Brown2001, horesh2006thesis, malone2014, goren2018}. EIT also has a variety of non-medical applications such as industrial process tomography and non-destructive testing, see e.g. \cite{TAPP200317,NISSINEN20141,Seppanen2001,Hallaji_2014}.

Several reconstruction methods for EIT exist, tracing back to the seminal works by Calder\'on \cite{Calderon1980} as well as Barber and Brown \cite{Barber1984}.  The majority of reconstruction methods can be classified into the following categories: linearized methods, Bayesian methods, 
optimization based methods, 
direct methods, and machine-learning based methods. We refer the reader to \cite{Holder2005, Mueller2012, Brown2009,Kaipio2000} for more thorough reviews of such methods.  Here we focus on complex-geometrical optics (CGO) solution based methods, a subset of direct methods for EIT, in particular Calder\'on and D-bar. Calder\'on's method ~\cite{Calderon1980} is a linearized CGO-based method that holds in both two and three dimensions.  It has been primarily studied from a theoretical standpoint, but has recently been shown to be suitable for medical imaging applications as well \cite{Bikowski2008, Boverman2009, Muller2017}.  The D-bar method, in contrast, is a fully non-linear CGO-based method developed for EIT in \cite{Novikov1987, Nachman1988a, Nachman1996}.  D-bar based reconstruction methods have been developed and implemented in 2D for over a decade (e.g., \cite{Isaacson2006,Knudsen2009,Dodd2014,Hamilton2013} with a recent push to directly embed a priori information into the reconstruction process \cite{AlsakerMueller2016,Alsaker2017} as well as handle partial boundary data \cite{Hauptmann2017,Hauptmann2017a}.  Progress in three-dimensions for D-bar methods has been slower to emerge and focused exclusively on theoretical studies and continuum-based numerical implementations \cite{Bikowski2011,Knudsen2011,Delbary2011,Delbary2014,Cornean2006}.  It is important to note that the EIT problem is inherently three-dimensional.  This is due to the fact that even if currents are injected into a single planar cross-section of a domain, the currents will have an out-of-plane effect inside the domain.  Therefore, transitioning to three-dimensional image reconstruction methods is of great interest.


To the best of our knowledge, this paper is the first to use a 3D D-bar method on physically realistic electrode data simulated with the Complete Electrode Model using the Finite Element method (FEM) \cite{Somersalo1992}.  Calder\'on's method has previously been used for electrode and experimental data before in a rectangular prism domain, \cite{Boverman2008a,Boverman2009} but this is the first time it has been used in a spherical domain with simulated CEM data and the computations are carried out in a different manner.

The manuscript is organized as follows: Section~\ref{sec:Methods} provides a general introduction to the EIT problem as well as overviews of the 3D D-bar method with $\texp$ and Calder\'on's method,  descriptions of the simulations/experiments performed, and evaluation metrics to be used.  In Section~\ref{sec:implement}, key aspects of the numerical implementations of the CGO-based reconstruction methods are presented, as well as the total variation (TV) and smoothness promoting (Smooth) regularized non-linear least squares (LS) methods to which the CGO-based reconstructions will be compared.  Results are shown in Section~\ref{sec:results} and a discussion provided in Section~\ref{sec:discussion}.  Conclusions are drawn in Section~\ref{sec:conclusion}.

\section{Methods}\label{sec:Methods}
The imaging task of recovering the internal admittivity of an object from electrical measurements taken at the surface uses the admittivity equation
\begin{equation}\label{eq:cond}
\nabla\cdot \gamma(x)\nabla u(x)  =  0, \quad x\in\Omega\subset\R^n,\quad n\geq 2,\\
\end{equation}
derived from Maxwell's equations, where $u(x)$ denotes the electrical potential and $\gamma(x)$ is the spatially dependent isotropic admittivity defined on a domain $\Omega\subset\R^n$ with Lipschitz boundary.  The admittivity $\gamma(x)=\sigma(x) + i\omega\varepsilon(x)\in L^\infty(\Omega)$ is comprised of the conductivity $\sigma$ and electrical permittivity $\varepsilon$ satisfying $\frac1C\leq \sigma\leq C$ for some constant $C$, where $\omega$ denotes the frequency of the applied current. In this work we consider the $n=3$ case.  We consider a Neumann boundary condition  $\gamma\frac{\partial u}{\partial\nu} = g(x)$, which corresponds to injecting current at the surface of the domain and measuring the resulting surface voltage $u(x)$ on $\bndry$.  Varying the current injection pattern and repeating surface voltage measurements corresponds to gaining approximate knowledge of the Neumann-to-Dirichlet (ND) map $\mathcal{R}_\gamma: \gamma\frac{\partial u}{\partial\nu} \mapsto u(x)$.  In practice, the CGO-based methods used in this paper rely on the Dirichlet-to-Neumann (DN) map, $\Lambda_\gamma=\mathcal{R}_\gamma^{-1}$. In this work, for simplicity, we assume $\gamma=\gammaB=1$ is constant in a neighborhood of the boundary.  See \cite{Nachman1988a} for the alterations required for the non-constant case.

We remark that the Neumann and Dirichlet boundary conditions used in the theoretical CGO-based methods constitute a simplified and approximate model for the physically realistic boundary conditions more often modeled by the complete electrode model~\cite{Somersalo1992} in optimization-based methods.  Further details regarding the discretiztaion of the DN and ND maps is discussed below in~\S\ref{sec-texp-implement}.

In this work we compare CGO-based reconstruction methods with well-established TV and smoothness promoting regularized non-linear least squares methods, highlighting that this is the first time that the $\texp$ D-bar method has been presented using physically realistic electrode data as well as the first direct comparison of CGO-based reconstruction methods and regularized non-linear least squares  methods in 3D.  We review the CGO-based methods studied below.

\subsection{The $\texp$ Method}\label{sec-texp}
The $\texp$ method is based on the fully non-linear D-bar reconstruction method introduced in \cite{Nachman1988a} and \cite{Novikov1987}.  The D-bar method for 3D EIT involves using a scattering transform, a sort-of non-linear Fourier transform, tailor-made for the Schr\"odinger equation.  This D-bar method was developed for the Schr\"odinger equation (see e.g., \cite{Beals1985}) and can be applied to EIT problem via the change of variables $\tilde{u} = \gamma^{1/2} u$ and $q(x) = \frac{\Delta \sqrt{\gamma(x)}}{\sqrt{\gamma(x)}}$,
\[\nabla\cdot \gamma(x)\nabla u(x)  =  0 \longrightarrow [-\Delta + q(x)] \tilde{u}(x) = 0.\]

Introducing an auxiliary variable, $\zeta\in \mathcal{V}_\xi$, special solutions, called {\it complex geometrical optics} solutions $\psi(x,\zeta)$ exist to 
\[ [-\Delta + q(x)] \psi(x,\zeta) = 0,\]
where $\psi(x,\zeta)\sim e^{ix\cdot\zeta}$ for large $|x|$ or $|\zeta|$.  We remark here that $\zeta$ is a purely nonphysical parameter.  The space $\mathcal{V}_\xi$ is a subspace of $\C^3$ defined by special orthogonality properties 
\begin{equation*}
\mathcal{V}_\xi = \left\{\zeta\in\C^3\middle|\zeta^2=0, \;\; \left(\xi+\zeta\right)^2=0\right\}, \quad \text{for each $\xi\in\R^3,$}
\end{equation*}
where $\zeta^2 = \zeta\cdot\zeta$.  The \emph{scattering transform}
\begin{equation}\label{eq:t-scat-full}
t(\xi,\zeta)=\int_{\R^3} e^{-ix\cdot(\xi+\zeta)}q(x)\psi(x,\zeta)\;dx,
\end{equation}
can be seen as the Fourier transform of the potential $q(x)$ if $|\zeta|$ is `large enough' by replacing $\psi(x,\zeta)$ with its asymptotic behavior $e^{ix\cdot\zeta}$
\begin{equation}\label{eq:t-scat-approx}
t(\xi,\zeta)\approx \int_{\R^3} e^{-ix\cdot(\xi+\zeta)}q(x)e^{ix\cdot\zeta}\;dx=\int_{\R^3}e^{-ix\cdot\xi} q(x)dx=\hat{q}(\xi),\quad \text{for $|\zeta|$ large.}
\end{equation}
Thus, if there were a way to determine the scattering data $t(\xi,\zeta)$ from the measured current/voltage data, one could then recover the conductivity by:
\[t \rightarrow \hat{q} \rightarrow q(x)  \rightarrow \gamma(x).\]
Alessandrini's identity \cite{Alessandrini1988} provides this missing link leading to
\begin{equation}\label{eq:t-scat-tBIE}
t(\xi,\zeta)=\int_{\bndry} e^{-ix\cdot(\xi+\zeta)}\left(\Lambda_\gamma-\Lambda_1\right)\psi(x,\zeta)\;dS(x),
\end{equation}
which instead requires knowledge of the CGOs $\psi(x,\zeta)$ for $x\in\bndry$ and $\zeta\in \mathcal{V}_\xi$.

The full non-linear D-bar method obtains those traces by solving an additional boundary integral equation involving a special Faddeev Green's function \cite{Faddeev1966} analogous to the 2D D-bar setting.  As this step is ill-posed and more computationally demanding, we proceed here with the `exp' approximation to the scattering data, denoted $\texp$, which replaces the traces of the CGO solutions $\psi$ with their asymptotic behavior $e^{ix\cdot\zeta}$ giving
\begin{equation*}
\texp(\xi,\zeta)=\int_{\bndry} e^{-ix\cdot(\xi+\zeta)}\left(\Lambda_\gamma-\Lambda_1\right)e^{ix\cdot\zeta}\;dS(x).
\end{equation*}
This can be considered a {\it `Born approximation'} to the fully non-linear scattering data.

The steps of the $\texp$ algorithm are then as follows:
\[\left(\Lambda_\gamma,\Lambda_1\right)\overset{1}{\longrightarrow}\texp(\xi,\zeta) \overset{2}{\longrightarrow}\qexp \overset{3}{\longrightarrow}\sigexp.\]

\noindent\framebox{\parbox{\textwidth}{\small
\begin{enumerate}[1.]
\item {\bf Compute the Approximate Scattering Data $\texp$.}  For each $\xi\in\R^3$, fix $\zeta\in\mathcal{V}_\xi$ and compute the scattering transform using the asymptotic behavior of the CGO solutions $\psi(x,\zeta) \sim e^{ix\cdot\zeta}$:
\begin{equation}\label{eq:texp-scat}
\texp(\xi,\zeta)=\int_{\bndry} e^{-ix\cdot(\xi+\zeta)}\left(\Lambda_\gamma-\Lambda_1\right)e^{ix\cdot\zeta}\;dS(x).
\end{equation}
where 
\begin{equation}\label{eq:Vxi}
\mathcal{V}_\xi = \left\{\zeta\in\C^3\middle|\zeta^2=0, \;\; \left(\xi+\zeta\right)^2=0\right\}, \quad \text{for each $\xi\in\R^3.$}
\end{equation}

\vspace{1em}

\item {\bf Compute the Approximate Potential $\qexp$:} Recover the approximate potential $\qexp$ from its Fourier transform $\widehat{\qexp}(\xi)$ by using the $\texp$ scattering data for large $|\zeta|$:  
\begin{equation}\label{eq:texp-q}
\qexp(x)= \mathcal{F}^{-1}\left\{\widehat{\qexp}(\xi)\right\}(x)\approx \mathcal{F}^{-1}\left\{\texp(\xi,\zeta)\right\}(x)=\frac{1}{(2\pi)^3}\int_{\R^3} e^{i x\cdot\xi} \texp(\xi,\zeta)\;d\xi, \quad x\in\R^3.
\end{equation}

\item {\bf Recover the approximate admittivity $\sigexp$:} Compute the approximate admittivity $\sigexp$ by solving the following boundary value problem:
\begin{equation}\label{eq:texp-QtoSigma}
\left\{\begin{array}{rclcl}
(-\Delta + \qexp(x)) \tilde{u}(x) & =& 0 & \qquad & x\in\Omega\subset\R^3\\
\tilde{u}(x) & =& 1 && x\in\bndry.
\end{array}
\right.
\end{equation}
Then $\sigexp(x)=\left(\tilde{u}(x)\right)^2$.
\end{enumerate}
}}

\normalsize
\vspace{1em}

Note that the algorithm as outlined above assumes that $\gamma(x)=\gammaB(x)\equiv 1$ near $\bndry$.  If this $\gammaB\neq1$, one can scale the DN map by instead using $\frac{1}{\gammaB}\Lambda_\gamma$, as in the 2D approach, and re-scaling at the end taking $\sigexp(x)=\gammaB\left(\tilde{u}(x)\right)^2$.  See \cite{Isaacson2004} for further details on scaling the DN map.
\subsection{The Calder\'on Method}\label{ssec-cald-orig}
Calder\'on's method, as considered here, comes from Calder\'on's original paper \cite{Calderon1980} and has a different flavor than the D-bar method above. Instead, it focuses on approximating the complex admittivity $\gamma(x)$ directly from its Fourier transform.  For discussion purposes we will call this the Calder\'on's $\Fhat$ ({\it F-hat}) method.  For further details the interested reader is referred to \cite{Calderon1980, Bikowski2008,Muller2017}. The main idea behind this method is to assume the complex admittivity is a constant plus a small perturbation, $\gamma(x)=\gammaB+\delta\gamma(x)$.

The steps of the algorithm are as follows:

\[\left(\Lambda_\gamma,\Lambda_b\right)\overset{1}{\longrightarrow}\Fhat(z) \overset{2}{\longrightarrow}\delta\sigcal\overset{3}{\longrightarrow}\sigcal.\]

\noindent\framebox{\parbox{\textwidth}{\small
\begin{enumerate}[1.]
\item {\bf Approximate the Fourier transform of the admittivity, $\hat{\gamma}(z)$:} The approximate Fourier transform of the small perturbation, $\delta\gamma(x)$, is given by
\begin{equation}\label{eq:CaldFhat}
\widehat{\delta\gamma}(z)\approx\Fhat(z):=-\frac{1}{2\pi^2|z|^2}\int_{\partial\Omega}e^{\pi i(z\cdot x)+\pi(a\cdot x)}\left(\Lambda_\gamma-\Lambda_b\right) e^{\pi i(z\cdot x)-\pi(a\cdot x)}dS(x),
\end{equation}
where $z$ and $a$ satisfy 
\begin{equation}\label{eq:za_conditions}
    z,a\in\R^3, |z|=|a|,\text{ and }z\cdot a=0.
\end{equation}
\begin{enumerate}[a.]
    \item {\bf (Optional) Average over choices of $a$:} For fixed $z$, the choice of $a$ to satisfy \eqref{eq:za_conditions} is not unique.  
    This can be leveraged as in \cite{Boverman2009} by defining a vector orthogonal to both $z$ some choice of $a$ satisfying \ref{eq:za_conditions}, denoted $a^\perp$, so that $a^\perp$ satisfies $z\cdot a^\perp=a\cdot a^\perp=0$ and $|a^\perp|=|z|=|a|$.  Then define

\begin{equation}\label{eq-U1U2}
U_1=e^{\pi i(z\cdot x)+\pi(\cos(\Theta)a\cdot x+\sin(\Theta)a^\perp\cdot x)}\text{ and }U_2=e^{\pi i(z\cdot x)-\pi(\cos(\Theta)a\cdot x+\sin(\Theta)a^\perp\cdot x)}
\end{equation}
 for $\Theta\in(0,2\pi]$ and average over $\Theta$. This leads to an alternate definition of $\Fhat(z)$:

\begin{equation}\label{eq:betterFhat}
\Fhat(z):=-\frac{1}{2\pi^2|z|^2}\frac{1}{2\pi}\int_0^{2\pi}\int_{\partial\Omega}U_1(z,x,\Theta)\left(\Lambda_\gamma-\Lambda_b\right) U_2(z,x,\Theta)dS(x)d\Theta.
\end{equation}

\end{enumerate}
\item {\bf Recover the perturbation in admittivity $\delta\sigcal$:} Compute the linearized perturbation in admittivity by taking the inverse Fourier transform of $\Fhat(z)$. 
\begin{equation}\label{eq:cal-FhattoGammaCart}
\delta\sigcal(x) \approx\mathcal{F}^{-1}\{\Fhat(z)\}(x)=\int_{\mathbb{R}^3}{\Fhat(z)e^{-2\pi i(x\cdot z)}dz}.
\end{equation}

\item {\bf Recover the admittivity, $\sigcal$:} The full admittivity can then be recovered by adding the background to the perturbation via
\begin{equation}\label{eq:sigCal_pert}
\sigcal(x)=\gammaB+\delta\sigcal(x).
\end{equation}

\end{enumerate}
}}

\normalsize
\vspace{1em}

We remark that while the D-bar method above uses a slightly different definition of the Fourier transform, both methods are presented using the definition consistent with their respective literature.





In this paper, we focus on spherical domains. On a sphere of radius $\rho$, let $x|_{\bndry}=\rho(\cos\phi\sin\theta,\sin\phi\sin\theta,\cos\theta)$. To satisfy \eqref{eq:za_conditions} we choose $z=|z|(\cos\tilde\phi\sin\tilde\theta,\sin\tilde\phi\sin\tilde\theta,\cos\tilde\theta)$ and $a=|z|(\cos\tilde\phi\cos\tilde\theta,\sin\tilde\phi\cos\tilde\theta,-\sin\tilde\theta)$ for $|z|\geq0$, $0\leq\tilde\phi\leq2\pi$ and $0\leq\tilde\theta\leq\pi$.  Converting \eqref{eq:cal-FhattoGammaCart} to spherical coordinates
yields
\begin{equation}\label{eq:cal-FhattoGamma}
\delta\sigcal(x) =\int_{0}^{\infty}{\int_0^{2\pi}{\int_0^\pi{|z|^2\sin\tilde\theta\Fhat(|z|,\tilde\phi,\tilde\theta)e^{-2\pi i(x\cdot z)}d\tilde\theta d\tilde\phi d|z|}}}.
\end{equation}

Additionally, Calder\'on~\cite{Calderon1980} proved that a mollifying function, $\hat{\eta}\left(\frac{z}{y}\right)$, for some $y\in\mathbb{R}$, can be applied via $\Fhat(z)\hat{\eta}\left(\frac{z}{y}\right)$ to reduce Gibbs phenomenon and maintain a good estimate for $\delta\sigcal$  

\begin{equation}\label{eq-sigcal-mollified}
\delta\sigcal(x) =\int_{0}^{\infty}{\int_0^{2\pi}{\int_0^\pi{|z|^2\sin\tilde\theta\Fhat(|z|,\tilde\phi,\tilde\theta)\hat{\eta\left(\frac{z}{y}\right)}e^{-2\pi i(x\cdot z)}d\tilde\theta d\tilde\phi d|z|}}}.
\end{equation}

We note that one can also directly reconstruct admittivity, $\sigcal(x)$ from a single data set by replacing $(\Lambda_\gamma-\Lambda_b)$ with $\Lambda_\gamma$ in the expression of $\Fhat$ and Step 2 will recover $\sigcal$ directly. To allow a more direct comparison of the D-bar and Calder\'on  reconstructions, the differential version \eqref{eq:sigCal_pert} will be used in this manuscript.

\subsubsection{The Calder\'on Method using $\texp$}\label{ssec-cald-texp}
We note that there is an alternative formulation of the Calder\'on method that uses the `Born' approximation to the scattering data $\texp$, described above in \S\ref{sec-texp}.  For discussion purposes, we will call this method the $\texpCAL$ method.  As shown in \cite{Bikowski2011}, this method involves a slight modification of the $\texp$ method replacing the two step process $\texp\rightarrow\qexp\rightarrow\sigexp$ with a single step $\texp\rightarrow\sigexpCal$, bypassing the potential $q$ completely.  The interested reader is referred to \cite{Bikowski2011} for further details.

The steps of the algorithm are as follows:
\[\left(\Lambda_\gamma,\Lambda_1\right)\overset{1}{\longrightarrow}\texp(\xi,\zeta) \overset{2}{\longrightarrow}\sigexpCal\]

\noindent\framebox{\parbox{\textwidth}{\small
\begin{enumerate}[1.]
\item {\bf Compute the Approximate Scattering Data $\texp$.}  For each $\xi\in\R^3$, fix $\zeta\in\mathcal{V}_\xi$ and compute the scattering transform $\texp$ via \eqref{eq:texp-scat}.

\vspace{1em}

\item {\bf Recover the approximate admittivity $\sigexpCal$:}  Compute the linearized Calder\'on admittivity $\sigexpCal$ by computing 
\begin{equation}\label{eq:cal-tEXPtoSigma}
\sigexpCal(x) = 1-\frac{2}{(2\pi)^3}\int_{\R^3}\frac{\texp(\xi,\zeta(\xi))}{|\xi|^2}e^{ix\cdot \xi} d\xi, \quad x\in\R^3, \quad \text{for $|\zeta|$ large.}
\end{equation}
\end{enumerate}
}}

\normalsize
\vspace{1em}
We note that this version of the Calder\'on method, and the D-bar method above in \S\ref{sec-texp}, require the difference in DN maps $\Lambda_\gamma - \Lambda_1$ and thus voltage data for $\gamma\equiv1$ must be obtained or simulated in some way to provide absolute EIT images.  However, as noted, Calder\'on's $\Fhat$ method can use data from a constant background, but it can  also produce absolute EIT images without any simulation of $\Lambda_b$.  Since we will show D-bar reconstructions which use $\texp$, we will focus on showing reconstructions from Calder\'on's $\Fhat$ method.  We note that if $\gammaB\neq1$, the DN map can be scaled to produce $\Lambda_1$ as in \cite{Isaacson2004} and described above in \S\ref{sec-texp}.

\subsection{Assessment of the Methods}\label{sec:methods-assess}
To investigate the quality of the $\texp$ and Calder\'on reconstruction methods on the CEM electrode data we tested each method on the several examples described below, and directly compared the reconstruction to two common regularized non-linear least squares methods, outlined in \S\ref{sec:refrecon}.

\subsubsection{Examples Considered}
\begin{figure}[t]
    \centering
    \includegraphics[width=320pt]{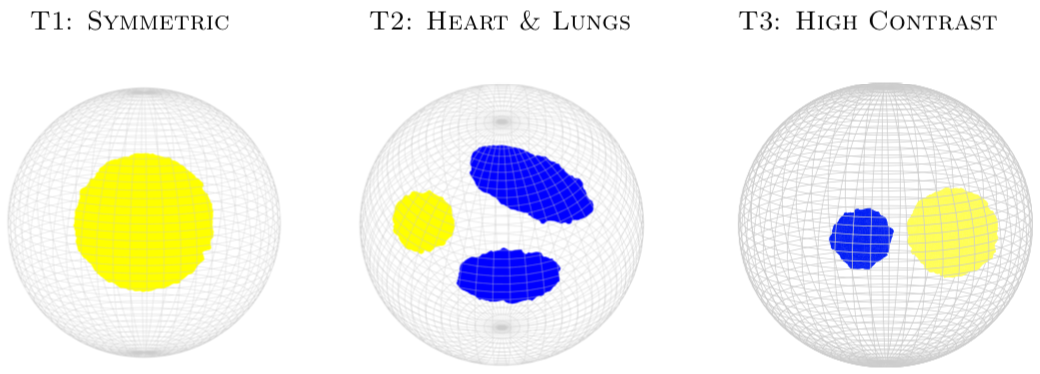}
%
%
%
%
%
    
    \caption{The simulated targets considered in this manuscript.}
    
    \label{fig:Examples}
\end{figure}

Figure~\ref{fig:Examples} displays the simulated targets used in this work.  We explore five questions using these targets: 
\begin{enumerate}
    \item How does reconstruction quality from electrode data compare to analytic data?
    \item Can admittivity targets be recovered?
    \item How is reconstruction quality affected by the number of electrodes simulated?
    \item Can high-contrast targets be recovered?
    \item What is the effect of noise (in the voltage data) on the reconstructed admittivity? 
\end{enumerate}

For each of these questions, the CGO (D-bar and Calder\'on) reconstructions are compared to reconstructions from the regularized non-linear least squares methods outlined in \S\ref{sec:refrecon}.  The domain used is a sphere of radius 1 meter for all examples.  Target T1 has a spherically symmetric ball of radius~0.5~m and conductivity~2~S/m contained in a unit sphere of conductivity $\gammaB=1$~S/m.  Target~T1 will be used to study Question~1.  Target T2 contains three disjoint targets crudely representing a heart and two lungs.  For this target we consider two different admittivity scenarios.  Target T2-A uses a heart with admittivity $2 + 0.6i$~S/m, and lungs $0.5 + 0.2i$~S/m with a nonunitary background admittivity of $0.8 + 0.3i$~S/m.  This complex-valued admittivity will be used to study the algorithms' capabilities with complex-valued data (Question~2).  As the heart and lungs target is commonly used in EIT algorithm literature, we also consider Target T2-B, a real-valued admittivity case with a conductive heart (2~S/m), resistive lungs (0.5~S/m) in a unitary background $\gammaB=1$~S/m.  Target T2-B will be used to study questions 3 and 5.  The high contrast target, T3, contains a conductive ball (1.5~S/m) of radius 0.3~m on the $x_1$ axis and smaller resistive ball (0.1~S/m) of radius 0.2~m on the $x_2$ axis with a background conductivity of $\gammaB=1$~S/m.  This target will be used to address questions 3-5.

The data was simulated using a finite element approximation of the complete electrode model (CEM) \cite{Cheng1989,Somersalo1992}. For details of the implementation see \cite{Vauhkonen1998a, Kaipio2000}. Following this approach, the domain $\Omega$ was divided into (507,767; 631,846; 632,710) tetrahedral elements and (93,880; 123,517; 124,775) nodes for the (32; 64; 128) electrode test cases, respectively. A measurement protocol with $L$ pairwise current injections where $L$ is the number of electrodes on the boundary $\partial \Omega$ was used, and for each injection the potentials on the electrodes were recorded, leading to $M =$~(1,024; 4,096; 16,384) simulated potential measurements. A current amplitude of 1~mA was used for the current injections and a contact impedance of 0.01~$\Omega$m$^2$ was used for all electrodes. To simulate noise in the data, relative mean Gaussian noise was added to the simulated noise-free data as follows:
\begin{equation}\label{eq:noisyVolts}
    V_{\mbox{\tiny noisy}}^j = V^j +  \eta\;  \texttt{mean}\left(|V^j|\right)N^j
\end{equation}
where $V^j$ denotes the voltage vector corresponding to the $j$th current pattern, $N^j$ a Guassian random vector unique for each current pattern $j$, and $\eta$ represents the noise level ($\eta=0.01$ corresponds to $1\%$ noise). Noise is not added to the data corresponding to a homogeneous background of $\sigma=1$ as that data would be purely simulated in experimental settings anyway for absolute EIT imaging with the D-bar and Calder\'on methods.


\subsubsection{Evaluation Metrics}\label{ssec:eval_metrics}
A quantitative assessment of Questions 1-5 will use several error metrics including localization and size of targets as well as whole-image error metrics.
Error metrics that are calculated per target are based on a segmentation using a threshold to identify targets, following \cite{Adler2009}.  For the complex case, T2-A, the real and imaginary parts are segmented separately. To segment, we first consider the difference from background of the image. Then we identify conductive targets as regions with greater than some threshold of the maximum difference from background conductivity and resistive targets as regions with less than some threshold of the minimum difference from background.  The conductive and resistive thresholds are set to 0.5 unless otherwise stated in \S\ref{ssec:eval_metrics_results}. Thresholds were only changed to better align the segmentation with how targets would be identified by eye from the reconstructions and occurred in cases with less signal (i.e. fewer electrodes and higher noise). For T2-A, the susceptivity thresholds are 0.8 and 0.3 respectively for the CGO methods and 0.5 and 0.3 for the LS-based methods due to the lower contrast between high and low susceptivity. A binary image is created in {\sc Matlab} from this segmentation, which is then labeled and regional properties such as the centroid location ($(x_1,x_2,x_3)$) and volume ($vol$) 
of each segmented target are computed using \texttt{regionprops3}.

Below are the error metrics we consider.
\begin{itemize}
\item \underline{\it Dynamic Range} (DR) - This is the ratio of the difference between the maximum and minimum values in the reconstructed image and that of the true image, 
    \begin{equation}\label{eq:DR}
        \text{DR}=\frac{\displaystyle\max_{x\in\Omega}\left\{\gamma_{recon}(x)\right\}-\min_{x\in\Omega}\left\{\gamma_{recon}(x)\right\}}{\displaystyle\max_{x\in\Omega}\left\{\gamma_{true}(x)\right\}-\min_{x\in\Omega}\left\{\gamma_{true}(x)\right\}}\times100\%.
    \end{equation}
    A perfect dynamic range is 100\%.
    
    \vspace{0.5em}
    \item \underline{\it Mean Square Error} (MSE) - This is computed using {\sc Matlab's} \texttt{immse} function within the whole spherical domain. The MSE is computed separately for the real and imaginary parts of the reconstruction. The MSE could be computed for the complex reconstruction, however, we want to focus on the quality of real and imaginary parts separately. A perfect MSE is 0.
    
    \vspace{0.5em}
    \item \underline{\it Multi-scale Structural Similarity Index} (MS-SSIM) - This is computed using {\sc Matlab's} \texttt{multissim3}, based on \cite{Wang2003}. Since the reconstructions of our spherical domain are stored within a rectangular prism, we fill values outside of the spherical domain with our estimate for the background conductivity, \eqref{eq-gam-best}, for the reconstruction and the true background conductivity, $\gammaB$, for the truth image. A perfect MS-SSIM is 1.
    
    \vspace{0.5em}
    \item \underline{\it Localization Error} (LE) - This is the distance between the reconstructed target's centroid, $\left(x_1^{recon},x_2^{recon},x_3^{recon}\right)$, and the true target's centroid, $\left(x_1^{true},x_2^{true},x_3^{true}\right)$,
    \begin{equation}\label{eq:LE}
        \text{LE}=\sqrt{(x_1^{recon}-x_1^{truth})^2+(x_2^{recon}-x_2^{truth})^2+(x_3^{recon}-x_3^{truth})^2}.
    \end{equation}
    LE$=0$ means the reconstructed target is located where the true target is. 
    
    \vspace{0.5em}
    \item \underline{\it Relative Volume Ratio} (RVR) - This is the ratio of reconstructed target's volume and the true target's volume. This is the same as what \cite{Adler2009} calls ``Amplitude Response," but renamed to reflect more clearly what is being measured. The relative volume ratio is computed by
    \begin{equation}\label{eq:RVR}
        \text{RVR}=\frac{vol_\text{reconstructed target}}{vol_\text{true target}},
    \end{equation}
    where $vol$ is given by \texttt{regionprops3} as the number of voxels in the segmented target. A perfect RVR is 1.
       
    \vspace{0.5em}
    
\end{itemize}
It is important to note that the LE and RVR evaluation metrics rely on this segmentation and thus are affected by the choice of threshold. For this reason, we have chosen to be as consistent as possible with thresholds for segmentation across all methods to compare their performance and any differences will be noted.

\subsubsection{Analytic Derivation of Continuum Data}\label{sec:AnalyticEval}
To address Question 1, we will compare Calder\'on reconstructions on analytic data given a radially symmetric conductive target in the sphere. This will then be compared to simulated electrode data of the same target (T1). Here we derive the analytic data and describe how Calder\'on's method can reconstruct the conductivity from that data. This will allow us to study the effect of moving from continuum boundary data to electrode data as well as the effect of decreasing the number of electrodes simulated (Questions 1 and 4).  Comparisons will be made in the frequency domain $\Fhat$ as well as for the reconstructed conductivity $\sigma(x)$.  

If the admittivity is real-valued and $\sigma(x)$ is radially symmetric, the analytic Fourier data $\Fhat$ can be computed via the eigenvalues of the DN maps $\Lambda_\sigma$ and $\Lambda_1$ as was done in \cite{Bikowski2011}.  
In a spherical domain with a radially symmetric, piece-wise constant, conductivity $\sigma$, the eigenfunctions of the DN map are the spherical harmonics.  
The corresponding eigenvalues
are given by the following proposition.

\begin{prop}[Proposition 3.6 of \cite{Bikowski2011}]\label{prop-3o6-BKM2011}
Suppose $\sigma$ is a radially symmetric piece-wise constant conductivity where the regions are defined radially by $0=r_0<r_1<r_2<\cdots<r_{N-1}<r_N=1$, and 
\begin{equation}\label{eq-sig-pw-rad}
\sigma(x)=\sigma_j>0, \qquad |x|\in[r_{j-1},r_j].
\end{equation}
Then, the eigenvalues of $\Lambda_\sigma$ are given by 
\begin{equation}\label{eq-evals-radsym-pw-constant}
\lambda_0 = 0,\qquad \lambda_\ell = \ell -\frac{2\ell+1}{1+C_{N-1}}, \qquad \ell>0.
\end{equation}
where $C_j=w_j\frac{\beta_\ell\sigma_{j+1}\rho_j + \sigma_j}{\sigma_{j+1}\rho_j - \sigma_j}$, with $\rho_1=1$, $\rho_j = \frac{C_{j-1} + w_j}{C_{j-1}-\beta_\ell w_j}$, $\beta_\ell = \frac{\ell+1}{\ell}$, and $w_j=r_j^{-(2\ell+1)}$.
\end{prop}
The eigenvalues of $\Lambda_1$ are given by $\lambda_{\ell}=\ell$.

Applying \cite{Bikowski2011} to  Calder\`on's method outlined in \S\ref{ssec-cald-orig}, $\Fhat$ can be computed using these eigenvalues via
\begin{equation}\label{eq-Fhat-DN-vals}
\Fhat(z)=-\frac{1}{2\pi^2|z|^2}\frac{1}{2\pi}\sum_{k=1}^{N_\Theta}\sum_{\ell=0}^\infty \sum_{m=-\ell}^{\ell} \tilde{a}_{\ell m}^\ast (z,\Theta_k) \tilde{b}_{\ell m}(z,\Theta_k) [\lambda_\ell-\ell],
\end{equation}
where $\frac{1}{2\pi}\sum_{k=1}^{N_\Theta}[\cdot]$ averages the boundary integral in \eqref{eq:betterFhat} over $N_\Theta$ vectors orthogonal to $z$. The coefficients $\tilde{a}_{\ell m}^\ast (z)$ and $\tilde{b}_{\ell m}(z)$ are defined by
\begin{eqnarray}
e^{\pi i(z\cdot x)+\pi(\cos(\Theta)a\cdot x+\sin(\Theta)a^\perp\cdot x)} &=& \sum_{\ell=0}^\infty \sum_{m=-\ell}^{\ell}  \tilde{a}_{\ell m}^\ast (z,\Theta) \left[Y^{m}_{\ell}(\theta,\phi)\right]^\ast \label{eq-coeffs-exp-zpa}\\
e^{\pi i(z\cdot x)-\pi(\cos(\Theta)a\cdot x+\sin(\Theta)a^\perp\cdot x)} &=& \sum_{\ell=0}^\infty \sum_{m=-\ell}^{\ell}  \tilde{b}_{\ell m}(z,\Theta) Y^{m}_{\ell}(\theta,\phi). \label{eq-coeffs-exp-zma}
\end{eqnarray}
Using the spherical harmonics $Y^{m}_{\ell}$ and their conjugates $\left(Y^m_\ell\right)^\ast$ allows the exploitation of their orthonormality reducing the boundary integral \eqref{eq:betterFhat} to \eqref{eq-Fhat-DN-vals} which only requires the coefficients of the exponential terms in the spherical harmonic basis. We then reconstruct the conductivity via equations and \eqref{eq:cal-FhattoGammaCart} and  \eqref{eq:sigCal_pert} above.

 It should be noted that $\lambda_\ell$ can be used in place of $[\lambda_\ell-\ell]$ in \eqref{eq-Fhat-DN-vals} to directly reconstruct the conductivity, but there is a significant Gibbs phenomenon at the domain boundary when this is done.   The expansions of the exponential terms in the spherical harmonics basis were computed using S2Kit to allow for fast computation of the coefficients of the spherical harmonic expansion.\footnote{The S2Kit package is hosted at \url{https://github.com/PatrickFBray/s2kit}.  We modified a Matlab MEX interface with the S2kit library to work with Matlab 2019a~\cite{Rodgers}.  The updated interface code is hosted at \url{https://github.com/PatrickFBray/s2kitmex}.}

\section{Numerical Implementation}\label{sec:implement}
In practice, $\gBest$, the best constant admittivity fit to the data, is often used since the true value of the admittivity in a neighborhood of the boundary $\bndry$ is unknown in practice.  This can be computed, solving a least-squares problem, as 
\begin{equation}\label{eq-gam-best}
\gBest=\frac{\sum_{k=1}^K\sum_{\ell=1}^LU_\ell^k(1)U_\ell^k(1)}{\sum_{k=1}^K\sum_{\ell=1}^LU_\ell^k(1)V_\ell^k},
\end{equation}
where $U_\ell^k(1)$ is the $k^{th}$ {\it simulated} voltage pattern measured on electrode $\ell$ with a homogeneous admittivity of $1$ and $V_\ell^k$ is the $k^{th}$ voltage pattern {\it measured} on electrode $\ell$ for the inhomogeneous admittivity $\gamma$~\cite{Isaacson2004}.  In an effort to simulate a more realistic experiment we treat $\gammaB$ as unknown and use the $\gBest$ approximation throughout.

Here we provide the numerical details pertinent to the implementation of the algorithms outlined above. 
\subsection{Implementation of the $\texp$ Method}\label{sec-texp-implement}
In this section we detail the numerical pieces needed to implement the $\texp$ algorithm described above in \S\ref{sec-texp}. 

Step~1 requires a discretized version of the integral in \eqref{eq:texp-scat}
\[\texp(\xi,\zeta)=\int_{\bndry} e^{-ix\cdot(\xi+\zeta)}\left[\left(\Lambda_\gamma-\Lambda_1\right)e^{i x\cdot\zeta}\right](x)\;dS(x).\]
Recall that this formulation requires that $\gamma\approx1$ near $\bndry$.  To ensure this is true, we computed the best constant admittivity fit to the measured data using $\gBest$ \eqref{eq-gam-best} and then considered the scaled DN map $\Lambda_{\tilde{\gamma}}$ corresponding to the scaled admittivity $\tilde{\gamma}\equiv\frac{\gamma}{\gBest}\approx 1$ near $\bndry$.  Following \cite{Isaacson2004}, $\Lambda_{\tilde{\gamma}}=\frac{1}{\gBest}\Lambda_{\gamma}$.  For simplicity of exposition, we drop the tilde notation and $\gamma$ will represent the scaled conductivity.  We can approximate the continuous integral for $\texp$ using a simple sum with quadrature points $x_\ell$ as follows
\begin{eqnarray}
\texp(\xi,\zeta)&=&\int_{\bndry} e^{-ix\cdot(\xi+\zeta)}\left[\left(\Lambda_\gamma-\Lambda_1\right)e^{i x\cdot\zeta}\right](x)\;dS(x)\nonumber\\
&\approx& \frac{4\pi}{L} \left[e^{-i\mathbf{x}\cdot(\xi + \zeta)}\right]^{\mbox{\tiny \bf T}}\mathbf{Q}\left(\mathbf{L}_\gamma - \mathbf{L}_1\right)\mathbf{Q}^{\mbox{\tiny \bf T}}\left[e^{i\mathbf{x}\cdot\zeta}\right],\nonumber
\end{eqnarray}
where $(\cdot)^{\mbox{\tiny \bf T}}$ denotes the traditional non-conjugate transpose,  $\mathbf{x}=[x_1,x_2,\ldots,x_L]^{\tiny \bf T}$ denotes the vector of $x$ values on the boundary $\bndry$ corresponding to the centers of the $L$ electrodes, $\mathbf{L}_\gamma$ and $\mathbf{L}_1$ denote the discrete matrix approximations to the DN maps $\Lambda_\gamma$ and $\Lambda_1$ respectively, $\mathbf{Q}$ denotes an orthonormal matrix created using the applied currents, and $\frac{4\pi}{L}$ is the uniform weight used for discretizing the surface area of the unit sphere.

The term $\mathbf{Q}\left(\mathbf{L}_\gamma - \mathbf{L}_1\right)\mathbf{Q}^{\mbox{\tiny \bf T}}\left[e^{i\mathbf{x}\cdot\zeta}\right]$ is used to approximate the action of the difference in DN maps $\Lambda_\gamma-\Lambda_1$ on the exponential behavior $e^{ix\cdot \zeta}$.  This has been done in the implementation of 2D D-bar methods by expanding the exponentials in the basis of applied current patterns, multiplying by the discrete matrix difference $\left(\mathbf{L}_\gamma - \mathbf{L}_1\right)$ and reforming the result by multiplying by the applied current pattern matrix.  Typically, in 2D this is done when applying trigonometric current patterns.  If pairwise current injection is used instead, then the `trig' voltage data is synthesized via a change of basis from the matrix of applied current patterns to the orthogonal trigonometric current patterns.  Note that the trigonometric current patterns are used as they are the eigenfunctions for the a circular domain in 2D.  The analog for the 3D case would be spherical harmonics.  


In this work we take a more general approach which avoids the transformation of the current patterns to the eigenfunction basis.  Instead, we produce an orthonormal matrix $\mathbf{Q}$ using a Modified Gram-Schmidt (MGS) algorithm on the set of linearly independent applied current patterns.  We note that there are several ways to produce a set of orthonormal vectors from a set of linearly independent vectors, but for this study MGS proved sufficient.  In this study pairwise current injection was simulated, in particular adjacent current patterns (skip-0) and thus $L-gcd(L,n_{\mbox{\tiny skips}}+1)=L-1$ linearly independent current patterns were applied where $gcd$ is the greatest common divisor and $n_{\mbox{\tiny skips}}$ denotes the number of electrodes skipped in the current injection scheme.  Letting $\mathbf{C}$ denote the $L\by(L-1)$ matrix of applied current patterns, the $L\by(L-1)$ matrix $\mathbf{Q}$ and upper triangular square $(L-1)\by(L-1)$ matrix $\mathbf{S}$ are produced using MGS, i.e. $\mathbf{C}=\mathbf{QS}$.  The measured voltage data is adjusted, if needed, so that each column sums to zero, i.e. each current injection produces voltages that sum to zero, and denoted $\mathbf{V}_{\gamma,{\mbox{\tiny meas}}}$.  Then the voltages that would have occurred if the orthonormal patterns $\mathbf{Q}$ were applied are approximated as
\begin{equation}\label{eq:v_synth}
    \mathbf{V}_{\gamma,{\mbox{\tiny synth}}} = \mathbf{V}_{\gamma,{\mbox{\tiny meas}}} \mathbf{S}^{-1}.
\end{equation}
Next, the discrete approximation $\mathbf{L}_\gamma$ to the continuous DN map $\Lambda_\gamma$ is formed as $\mathbf{L}_\gamma=\left(\mathbf{R}_\gamma\right)^{-1}$ where
\begin{equation}\label{eq:ND_discrete}
    \mathbf{R}_\gamma=\left[V_{\gamma,{\mbox{\tiny synth}}} \right]^*\mathbf{Q},
\end{equation}
using the inner product definition of the ND map, where $(\cdot)^*$ denotes the conjugate transpose.

To form the action $\Lambda_1 e^{i x\cdot\zeta}$ we used simulated voltage data corresponding to an admittivity of $\gamma=1$ computed using FEM with the Complete Electrode Model \cite{Somersalo1992, Kaipio2000}, producing $\mathbf{L}_1$ in the same manner as $\mathbf{L}_\gamma$.

The next important piece in computing \eqref{eq:texp-scat} is the auxiliary variable $\zeta$.  Following \cite{Delbary2014} we note that $\zeta\in\mathcal{V}_\xi$ forces $\zeta$ to have the following form:
\begin{equation}\label{eq:zeta-xi}
\zeta(\xi)=\left(-\frac{\xi}{2}+\left(\kappa^2-\frac{|\xi|^2}{4}\right)^{1/2}\xi^\perp\right)+i\kappa \xi^{\perp\perp},
\end{equation}
with $\kappa\geq\frac{|\xi|}{2}$ and $\xi^\perp,\;\xi^{\perp\perp}\in\R^3$ orthonormal vectors that are orthogonal to $\xi$.  Here we used $\kappa=\frac{|\xi|}{2}$, the minimal-zeta approach outlined in \cite{Delbary2014} that satisfies $\kappa\geq\frac{|\xi|}{2}$.  The orthonormal vectors $\xi^\perp,\;\xi^{\perp\perp}$ were computed for each corresponding $\xi$ using the the {\tt null} command in {\sc Matlab}.    Then, $\texp$ is computed on a finite $\xi$ grid $[-T_\xi,T_\xi]^3$ via
\begin{equation}\label{eq:texp-clipped}
    \texp(\xi,\zeta(\xi)) = \begin{cases}
    \frac{4\pi}{L} \left[e^{-i\mathbf{x}\cdot(\xi + \zeta)}\right]^{\mbox{\tiny \bf T}}\mathbf{Q}\left(\mathbf{L}_\gamma - \mathbf{L}_1\right)\mathbf{Q}^{\mbox{\tiny \bf T}}\left[e^{i\mathbf{x}\cdot\zeta}\right] & |\xi|\leq T_\xi\\
    0 & \text{else},
    \end{cases}
\end{equation}
with $\zeta(\xi)$ computed from \eqref{eq:zeta-xi}.

Next, in Step~2 we computed $\qexp$ from \eqref{eq:texp-q} as
\[\qexp(x)=\frac{1}{(2\pi)^3}\int_{[-T_\xi,T_\xi]^3} e^{i x\cdot\xi} \texp(\xi,\zeta(\xi))\;d\xi,\]
using a 3D Simpson's rule\footnote{A modified version of \cite{Padden2008} was used.}, noting that the integral over $\R^3$ reduces to the $\xi$-grid $[-T_\xi,T_\xi]^3$ as $\texp(\xi,\zeta)=0$ for $|\xi|>T_\xi$.

The boundary value problem \eqref{eq:texp-QtoSigma} in Step~3 was solved using the {\sc PDE toolbox} in {\sc Matlab} with the {\tt PDEmodel} structure using $\qexp$ with $14,014$ 3D  quadratic elements.  The approximate conductivity was recovered as $\sigexp(x)=\gBest\left(\tilde{u}(x)\right)^2$ where $\tilde{u}$ is the solution using the {\sc PDE toolbox}.  The solution was interpolated to a $128\by128\by128$ rectangular grid for viewing purposes using {\tt scatteredInterpolant}.

\subsection{Implementation of Calder\'on's Method}\label{sec-cal-Fhat-implement}
In this section, we detail the numerical pieces required to implement Calder\'on's method. As stated in \S\ref{ssec-cald-texp}, there is a connection between $\texp$ and Calder\'on's method.  Therefore, some of the implementation details from \S\ref{sec-texp-implement} will be the same or similar. For brevity, we will refer to the corresponding details from that section when necessary.

For Step~1, we will use the \eqref{eq:betterFhat} formulation of $\Fhat$. We discretize the boundary integral, denoted $I(z,\Theta)$ below, by

\begin{eqnarray}
I(z,\Theta)&=&\int_{\partial\Omega}U_1(z,x,\Theta)(\Lambda_\gamma-\Lambda_1) U_2(z,x,\Theta)dS(x)\nonumber\\
&\approx& \frac{\Delta\phi\Delta\theta}{A_e}U_1(z,\mathbf{x},\Theta)^{\mbox{\tiny \bf T}}\mathbf{Q}\left(\mathbf{L}_\gamma-\mathbf{L}_1\right)\mathbf{Q}^{\mbox{\tiny \bf T}}\left[U_2(z,\mathbf{x},\Theta)\right],\label{eq-discrete-Fhat_bdry}
\end{eqnarray}
where $\Delta\phi$ and $\Delta\theta$ are the minimum angles between electrode centers in the polar and azimuthal coordinates, respectively; $A_e$ is the uniform area of an electrode; $\mathbf{x}$ is a $27\times27\times27$ equally-spaced rectangular grid in this work; $\mathbf{Q}$ and $\mathbf{L}_\gamma$ are the same as in \S\ref{sec-texp-implement}; and $U_1$ and $U_2$ are given by \eqref{eq-U1U2}. We compute \eqref{eq-discrete-Fhat_bdry} at $\Theta_k\in(0,2\pi]$ for $k=1,\cdots,N_\Theta$, with $N_\Theta=30$ in all reconstructions in this paper.  We then averaged the computation over $\Theta$, in the same manner as in \eqref{eq-Fhat-DN-vals},

\begin{equation}\label{eq-Fhat-data}
\Fhat(z)\approx-\frac{1}{2\pi^2|z|^2}\frac{1}{2\pi}\sum_{k=1}^{N_\Theta}I(z,\Theta_k).
\end{equation}

Next, we computed the mollifier defined in ~\cite{Bikowski2008},
\begin{equation}\label{eq-Cald-moll}
\hat{\eta}\left(\frac{z}{y}\right)=e^{-\pi t|z|^2},
\end{equation}
where $y=1/\sqrt{t}$ and $t$ is treated as a regularization parameter. Setting $t=0$ implies there is no mollifying effect, and larger $t$ values will smooth the reconstruction and reduce large jumps at points of discontinuity in the admittivity.

Another regularization parameter arises in the truncation of the radial coordinate of the Fourier domain in the computation of \eqref{eq:betterFhat} to reduce the effect of noise in the data.  So, instead of letting $0\leq|z|<\infty$, we let $0\leq|z|<T_z$. This essentially has the same effect as the $\texp$ truncation in \eqref{eq:texp-clipped}. Thus, we reconstruct the perturbation from background conductivity by 

\begin{equation}\label{eq-cal-FhattoGamma-reg}
\delta\sigcal(x) =\int_{0}^{T_z}{\int_0^{2\pi}{\int_0^\pi{|z|^2\sin\tilde\theta\Fhat(|z|,\tilde\phi,\tilde\theta)e^{-\pi t|z|^2}e^{-2\pi i(x\cdot z)}d\tilde\theta d\tilde\phi d|z|}}}.
\end{equation}
The integration is computed using a 3D Simpson's rule with the $|z|, \tilde{\theta}$, and $\tilde{\phi}$ grids with $N_{|z|}, N_{\tilde{\theta}}$, and $N_{\tilde{\phi}}$ nodes on each axis, respectively.  For all real targets (T1, T2-B, T3), we chose $N_{|z|}=10$, $N_{\tilde{\theta}}=8$, and $N_{\tilde{\phi}}=14$.  For the complex target (T2-A), we chose $N_{|z|}=10$, $N_{\tilde{\theta}}=10$, and $N_{\tilde{\phi}}=20$.  Increasing the number of nodes in the Fourier domain can reduce artefacts, but also increases computation time. 

The reconstructions of $\sigcal(x)$ in this paper are then produced using \eqref{eq:sigCal_pert} replacing $\gamma_b$ with $\gBest$
\begin{equation*}
\sigcal(x)=\gBest+\delta\sigcal(x),
\end{equation*}
where $\gBest$ is given by \eqref{eq-gam-best}. The solution is then interpolated to a $128\times128\times128$ rectangular grid. 

\subsection{Regularized Non-Linear Least Squares}\label{sec:refrecon}
To compare the D-bar method to more common numerical 
3D absolute reconstructions, we include
reconstructions using a regularized non-linear Least Squares (LS) formulation. 
As many widely used regularization functionals and optimization techniques are based on real-valued variables, a common approach for complex-valued problems is to split the complex variables into real and imaginary parts. Utilizing such a real-valued formulation for the complex-valued EIT problem, the
(discretized) regularized LS approach amounts to finding the solution
\begin{equation}
  (\hat\sigma,\hat \epsilon)= {\rm arg} \min_{\sigma , \epsilon >0}\{\Vert V_{{\rm s}}-U_{{\rm s}}(\sigma,\epsilon)\Vert^2
+ \Psi (\sigma, \epsilon)\},
\label{lssol}
\end{equation}
where $\sigma \in \R^N$ is finite dimensional approximation of the conductivity and $\epsilon = \omega \varepsilon \in \R^N$ the susceptivity, leading to admittivity $\gamma = \sigma + i \epsilon \in \mathbb{C}^N$. The measurement vector $V_{{\rm s}} = ({\rm re}(V),{\rm im}(V))^{T}$ contains the real and imaginary parts of the measured voltages, 
$U_{{\rm s}}(\sigma,\epsilon)=({\rm re}(U(\sigma,\epsilon)),{\rm im}(U(\sigma,\epsilon))^T$ 
contains the real and imaginary parts of the
complex-valued forward model $U(\sigma,\epsilon) \in \mathbb{C}^{M}$ 
and $\Psi (\sigma, \epsilon)$ is a regularization functional. 

The implementation of the forward model $U(\sigma, \epsilon)$ is based on the Finite Element method (FEM) approximation of the complete electrode model (CEM) \cite{Somersalo1992}. In the FEM approximation, the electrical potential $u(x)$  
and the admittivity $\gamma (x)$ were approximated in separate tetrahedral FE meshes using piece-wise linear basis elements. The mesh for the admittivity is a nearly uniform mesh with $N = 36,364$ nodes and 195,948 elements, and the unknown admittivity was approximated in the piece-wise linear basis as $\gamma (x) = \sum_{k=1}^N (\sigma_k + i \epsilon_k) \psi_k (x)$, leading to the parameterizations $\sigma \in \R^N$ and $\epsilon \in \R^N$. The mesh for the electric potential is significantly more dense and refined near the electrodes with total of 94,980 nodes and 481,975 elements to achieve sufficient accuracy for the forward solution $U(\sigma,\epsilon)$. 

In this work, two different regularization functionals were considered.  The first is a smoothness promoting regularization 
\begin{equation}\label{l2psi}
\Psi (\sigma,\epsilon) = \| L_{\sigma} \sigma \|^2 + \| L_{\epsilon} \epsilon \|^2,
\end{equation}
where the regularization matrices are defined by $L_{\sigma}^{\mathrm{T}}L_{\sigma} = \Gamma_{\sigma}^{-1}$ and $L_{\epsilon}^{\mathrm{T}}L_{\epsilon} = \Gamma_{\epsilon}^{-1}$, and $\Gamma_{\sigma}$ and $\Gamma_{\epsilon}$ are smoothness promoting covariance matrices constructed using the correlation model \cite{Lieberman10}
\begin{equation}
\label{equ:GammaSmoothness}
  \Gamma(i,j) = a\exp\left(-\frac{\|x_i - x_j \|^2}{2b^2}\right) + c\delta_{ij}, \ \ \ i, j = 1,\ldots, N.
\end{equation}
where the correlation is based on the Euclidean distance between the positions $x_i$ and $x_j$. The positive scalar variables $a$, $b$ and $c$ are such that $a$ defines the prior variance of $\sigma$ or $\epsilon$, $b$ controls the correlation length and $c$ together with the Kronecker delta $\delta_{ij}$ guarantees that the covariance matrices are positive definite, allowing the computation of the Cholesky factors of the respective precision matrices. Parameter $b$ can be tuned by setting the distance $\|x_i - x_j \|$ to a desired value, e.g. half the radius of the target, and setting $\Gamma(i,j)$ to the desired covariance for that distance, e.g. 1\% of variance.

The second regularization model considered is a (smooth) total variation (TV), which is a well-established regularization model for
promoting blocky images thus appropriate for 3D admittivities where materials with different admittivities are distributed in well-defined objects with crisp boundaries \cite{Rudin1992}. Utilizing the piece-wise linear FE discretization, the smooth TV functional becomes
\begin{equation}\label{stvpsi}
\Psi (\sigma,\epsilon) = \alpha_{\sigma} \int_{\Omega} 
\sqrt{\| \nabla \sigma \|^2 +  \beta_{\sigma} }\ {\rm d} x
+
\alpha_{\epsilon} \int_{\Omega} 
\sqrt{\| \nabla \epsilon \|^2 +  \beta_{\epsilon} }\ {\rm d} x 
\end{equation}
where $\nabla \sigma$ and $\nabla \epsilon$ are the gradients of the conductivity and the susceptivity,
$\alpha_{\sigma} > 0$ and $\alpha_{\epsilon} > 0$ are regularization parameters and $\beta_{\sigma}$ and $\beta_{\epsilon}$ are smoothing parameters used to make the functional globally differentiable.

The non-linear optimization \eqref{lssol} was solved by a lagged Gauss-Newton method \cite{Vogel2002} equipped with a line search algorithm. The line search is implemented using bounded minimization such that the non-negativity $\sigma > 0$ and $\epsilon > 0$ is enforced using the log-barrier method \cite{Nocedal2006}. In the TV regularization, the regularization parameters $\alpha_{\sigma}$ and $\alpha_{\epsilon}$ were tuned manually for the best visual quality of the reconstruction.


\section{Results}\label{sec:results}
Here, reconstructions of the example admittivities and assessment metrics are presented for the driving questions outlined in \S\ref{sec:methods-assess}.  The regularization parameters used for the methods are stated below for each example.  All reconstructions were interpolated to a $128\by128\by128$ $x$-grid for uniformity in plotting.

\subsection{Electrode data vs. Analytic Data}
\begin{figure}[h!]
\centering
\includegraphics[width=335pt]{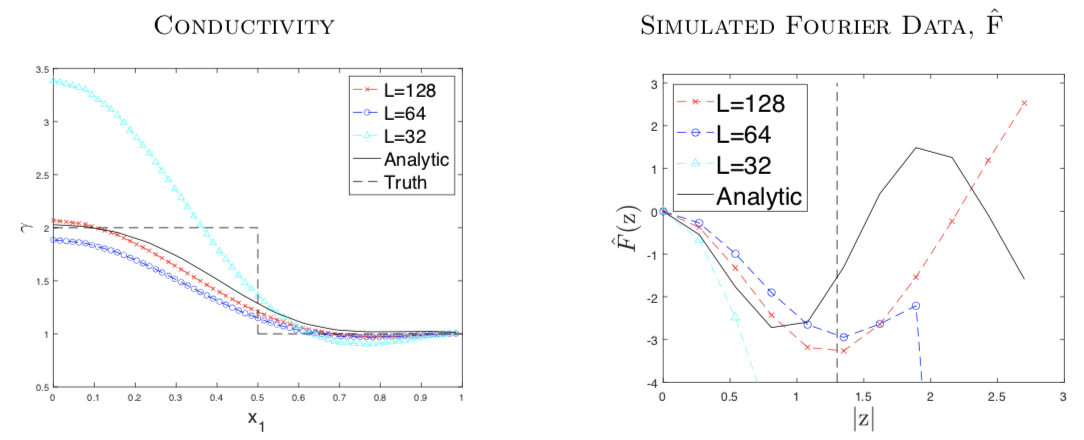}
%
%
%
%
%
\caption{\label{fig:T7-compare-rad} Comparison of reconstructed conductivity (Left) and Fourier data, $\Fhat$, (Right) for T1 using Calder\'on's method (equation \eqref{eq:sigCal_pert}). $T_z=2.7$ for the analytic data and $T_z=1.3$ for all three simulated electrode data cases. The mollifying parameter is $t=0.1$ for both analytic and simulated electrode data. The vertical dashed line indicates where the Fourier domain was truncated for the simulated electrode data cases.} 

\end{figure}

\begin{figure}[h!]
\centering
\begin{picture}(500,195)
\put(-5,0){\includegraphics[width=490pt]{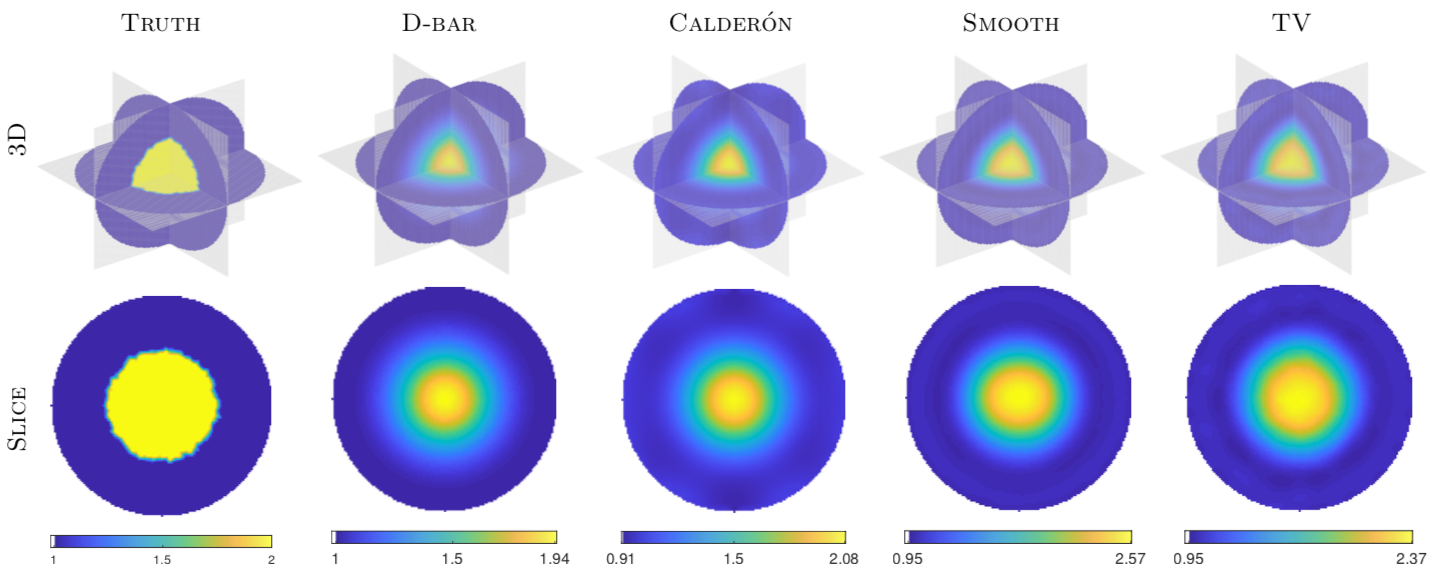}}
\end{picture}
\caption{\label{fig:T7-compare} Reconstructions of radially symmetric example T1 across algorithms using $L=128$ electrodes shown in 3D and a representative $x_2x_3$ slice.}
\end{figure}

Here we examine the disjoint target T1 from Fig.~\ref{fig:Examples}.  The Calder\'on reconstruction for analytic data was computed according to \S\ref{sec:AnalyticEval}, where we let the sums over $\ell$ stop at $\ell=50$ in equations \eqref{eq-Fhat-DN-vals}, \eqref{eq-coeffs-exp-zpa}, and \eqref{eq-coeffs-exp-zma}. Equation \eqref{eq-cal-FhattoGamma-reg} was used to reconstruct the conductivity's perturbation from background with Fourier regularization parameter $T_z=2.7$, and mollifying regularization parameter $t=0.1$. In addition, we reconstructed simulated 128, 64, and 32 electrode data using the implementation described in \S\ref{sec-cal-Fhat-implement} with Fourier regularization parameter $T_z=1.3$ and mollifying regularization parameter $t=0.1$ across all electrode configurations.  Since reconstructions from radially symmetric targets are radially symmetric, we show these Calder\'on reconstructions along the $x_1$-axis on the left of Figure~\ref{fig:T7-compare-rad}.  On the right of Figure~\ref{fig:T7-compare-rad}, we show a radial cross-section of $\Fhat$ across all data types, where a vertical dashed line highlights the Fourier regularization parameter $T_z$ used for the simulated electrode data reconstructions shown on the left.  Additionally, reconstructions in 3D planar cross-sections along with the $x_2x_3$ plane are shown in Figure~\ref{fig:T7-compare} for $\texp$, Calder\'on, and the two LS-based methods for the $L=128$ electrode case. The Calder\'on reconstructions in Figure~\ref{fig:T7-compare} correspond to the red dashed line with $\times$ markers for the $L=128$ case in Fig.~\ref{fig:T7-compare-rad}.

Reconstruction parameters for Fig.~\ref{fig:T7-compare} were as follows: $\texp$ ($T_\xi =7$), $\Fhat$ ($T_z=1.3, t=0.1$), Smooth ($a = 0.444,\; b = 0.330,\; c = 10^{-8}$), TV ($\alpha_{\sigma} = 2000$,\; $\beta_{\sigma} = 10^{-8}$).

\clearpage


\subsection{Complex-Valued Admittivity}
\begin{figure}[h!]
\centering
\includegraphics[width=425pt]{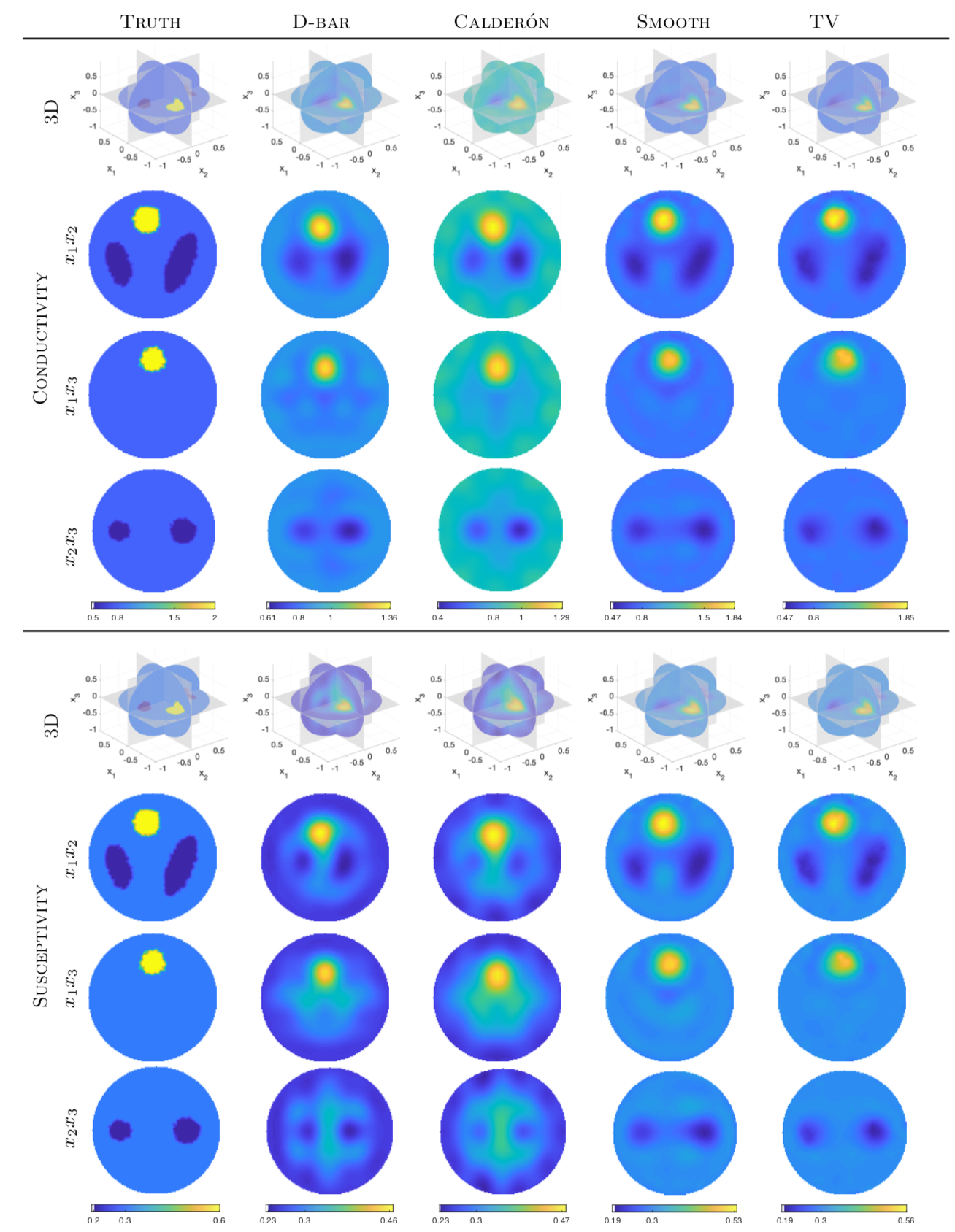}
\caption{\label{fig:T3cmplx-compare} Comparison of conductivity and susceptivity reconstructions for the complex-valued heart and lungs target T2-A. }
\end{figure}

We proceed with a common target for EIT reconstruction, the simulated heart and lungs phantom T2 in Figure~\ref{fig:Examples} using a complex-valued admittivity.  Reconstructions are shown in Figure~\ref{fig:T3cmplx-compare}.  The regularization parameters used were: $\texp$ ($T_\xi =14$), \newline $\Fhat$ ($T_z=2.3, t=0.05$), Smooth ($a_{\sigma} = 0.284,\; b_{\sigma} = 0.330,\; c_{\sigma} = 10^{-8};\; a_{\epsilon} = 0.040,\; b_{\epsilon} = 0.330,\; c_{\epsilon} = 10^{-8}$), TV ($\alpha_{\sigma} = 20,\; \beta_{\sigma} = 10^{-8},\; \alpha_{\epsilon} = 10,\; \beta_{\epsilon} = 10^{-8}$).

\subsection{Varying the Number of Electrodes}
\begin{figure}[h!]
\centering
\includegraphics[width=425pt]{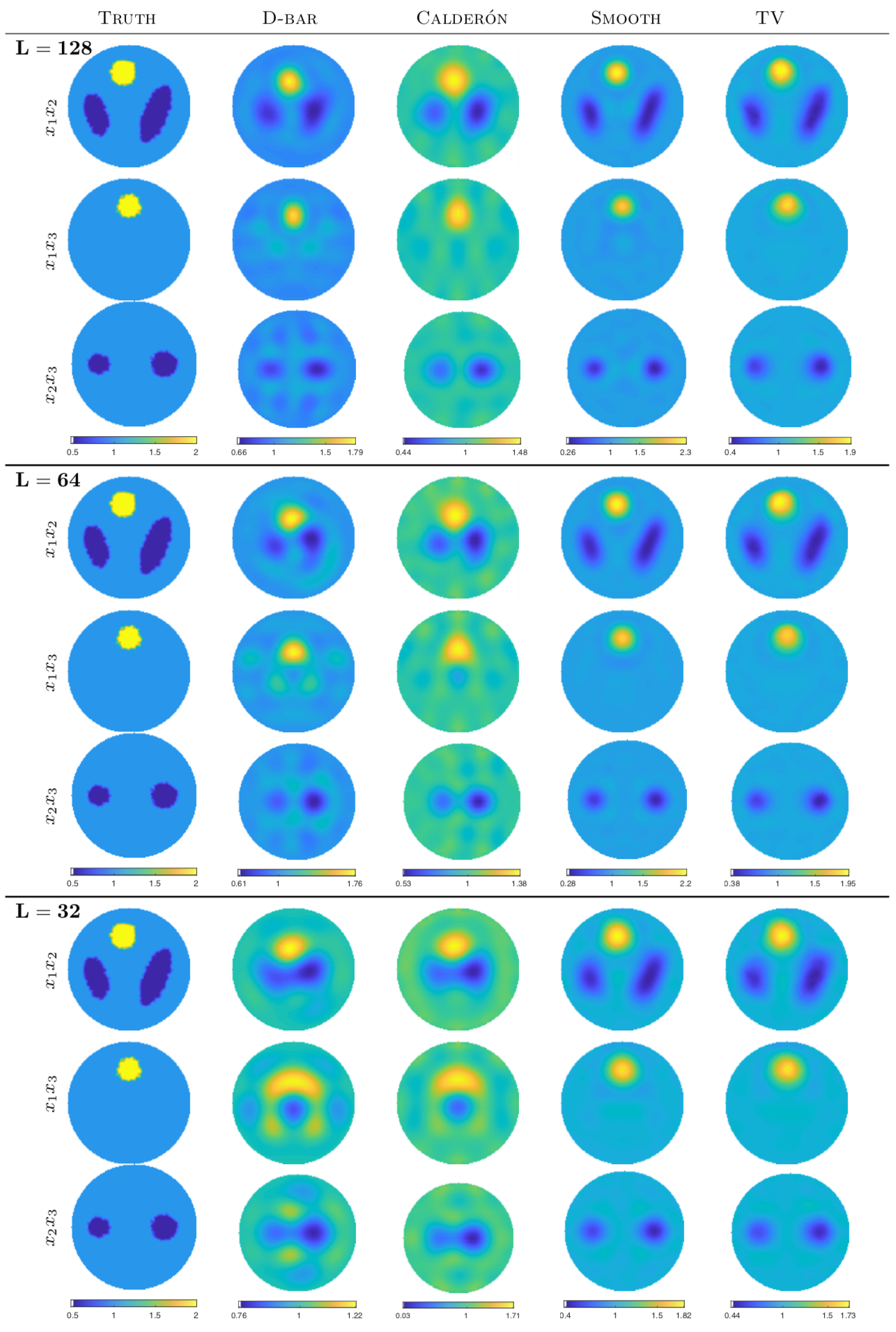}
\caption{\label{fig:T3-compare} Comparison of reconstructions for the real-valued heart and lungs target T2-B using $L=128$, $64$, or $32$ electrodes.}
\end{figure}

We now examine how the reconstructions change with more or fewer electrodes simulated on the boundary.  We begin by investigating the effect on the the real-valued heart and lungs target T2-B. Figure~\ref{fig:T3-compare} shows reconstructions using the reconstruction methods for the real-valued target T2-B from Fig.~\ref{fig:Examples} while varying the number of electrodes used: $L=128$, $L=64$ and $L=32$.   Each electrode has radius 0.05m corresponding to electrodes covering 8\%, 4\% and 2\% of the domain surface, respectively.  The regularization parameters were: $\texp$ ($T_\xi =16,\;14,\; 11$), $\Fhat$ ($T_z=2.3,\;2.3, 1.8$, $t=0.05$), Smooth ($a = 0.444,\;0.444,\;0.444; \; b = 0.165,\;0.165,\;0.165; \; c = 10^{-8},\;10^{-8},\;10^{-8}$), TV ($\alpha_{\sigma} = 200,\; 100,\; 50$, \; $\beta_{\sigma} = 10^{-8}$) for $L=128,\;64,$ and 32 electrodes, respectively.

\subsection{High-Contrast Targets}
\begin{figure}[h!]
\centering
\includegraphics[width=425pt]{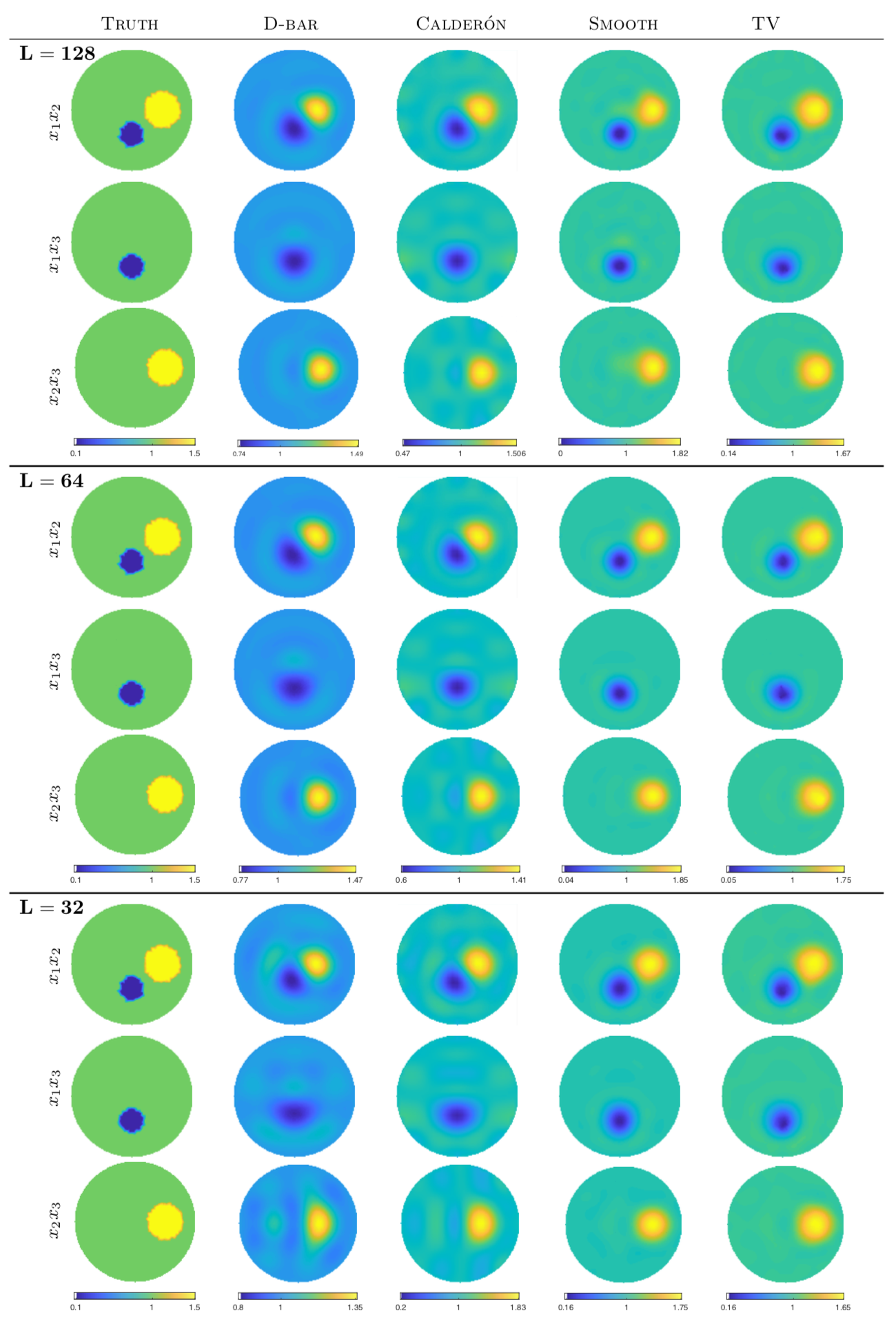}
\caption{\label{fig:T6-compare} Comparison of reconstructions for the high contrast target T3 using $L=128$, $64$, or $32$ electrodes.}
\end{figure}

We next explore the effect of adjusting the number of electrodes with a high contrast admittivity where the conductor ($\sigma=1.5$) is {\it fifteen times} as conductive as the resistor ($\sigma=0.1$).  Figure~\ref{fig:T6-compare} shows the reconstructions with  $L=128$, $64$, and $32$ electrodes on target T3.  The regularization parameters were: $\texp$ ($T_\xi =12,\;11, 11$), $\Fhat$ ($T_z=2,\;2, 1.8$, $t=0.05$), Smooth ($a = 0.444,\;0.444,\;0.444; \; b = 0.165,\;0.165,\;0.165; \; c = 10^{-8},\;10^{-8},\;10^{-8}$), TV ($\alpha_{\sigma} = 200,\; 100,\; 50$, \; $\beta_{\sigma} = 10^{-8}$) for $L=128,\;64,$ and 32 electrodes, respectively.  

\subsection{Noisy Voltage Data}
\begin{figure}[h!]
\centering
\includegraphics[width=425pt]{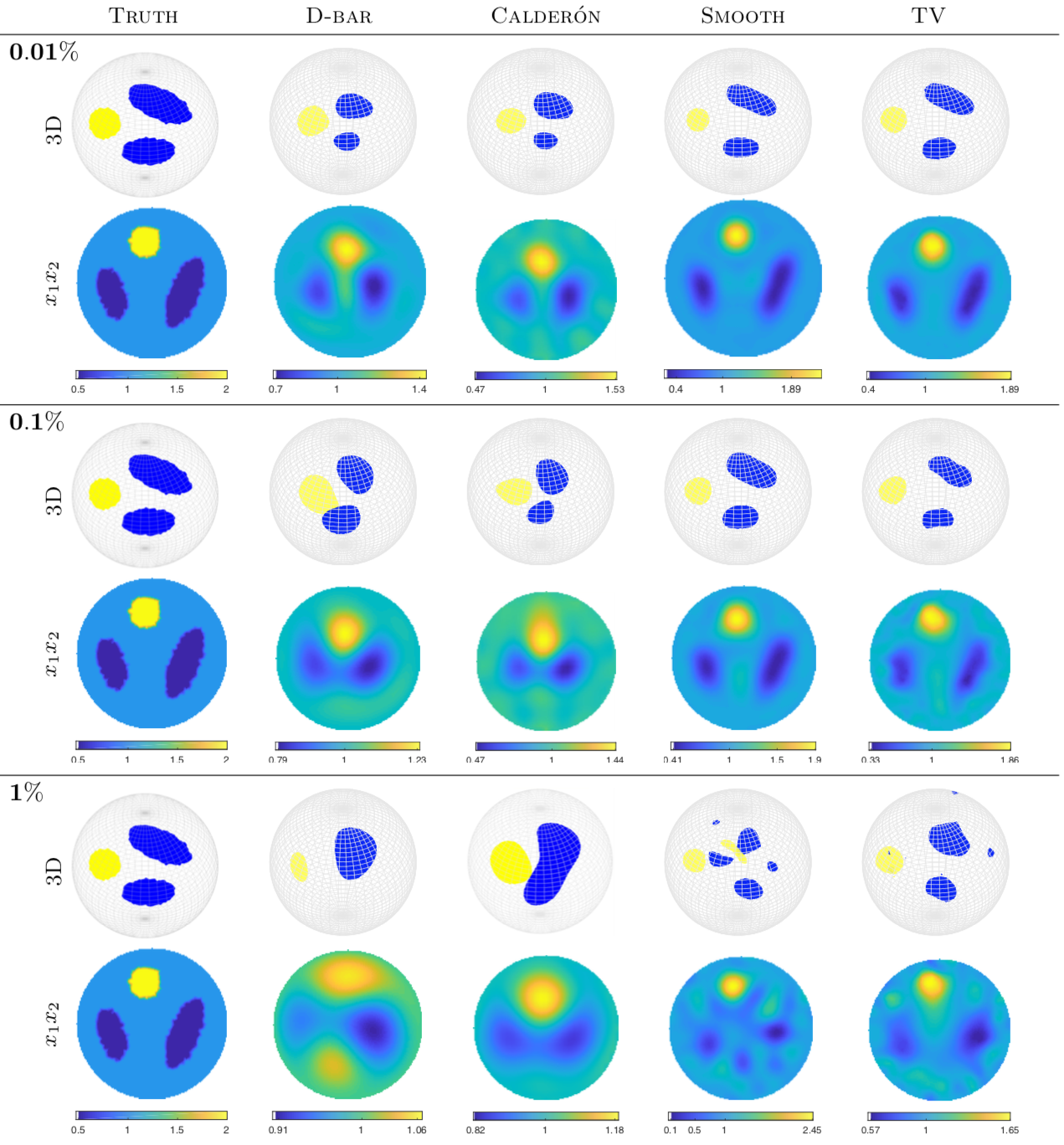}
\caption{\label{fig:T3_noisy} Comparison of reconstructions for the real-valued heart and lungs target T2-B with increasing levels of noise added to the voltage data.}
\end{figure}

\begin{figure}[h!]
\centering
\includegraphics[width=425pt]{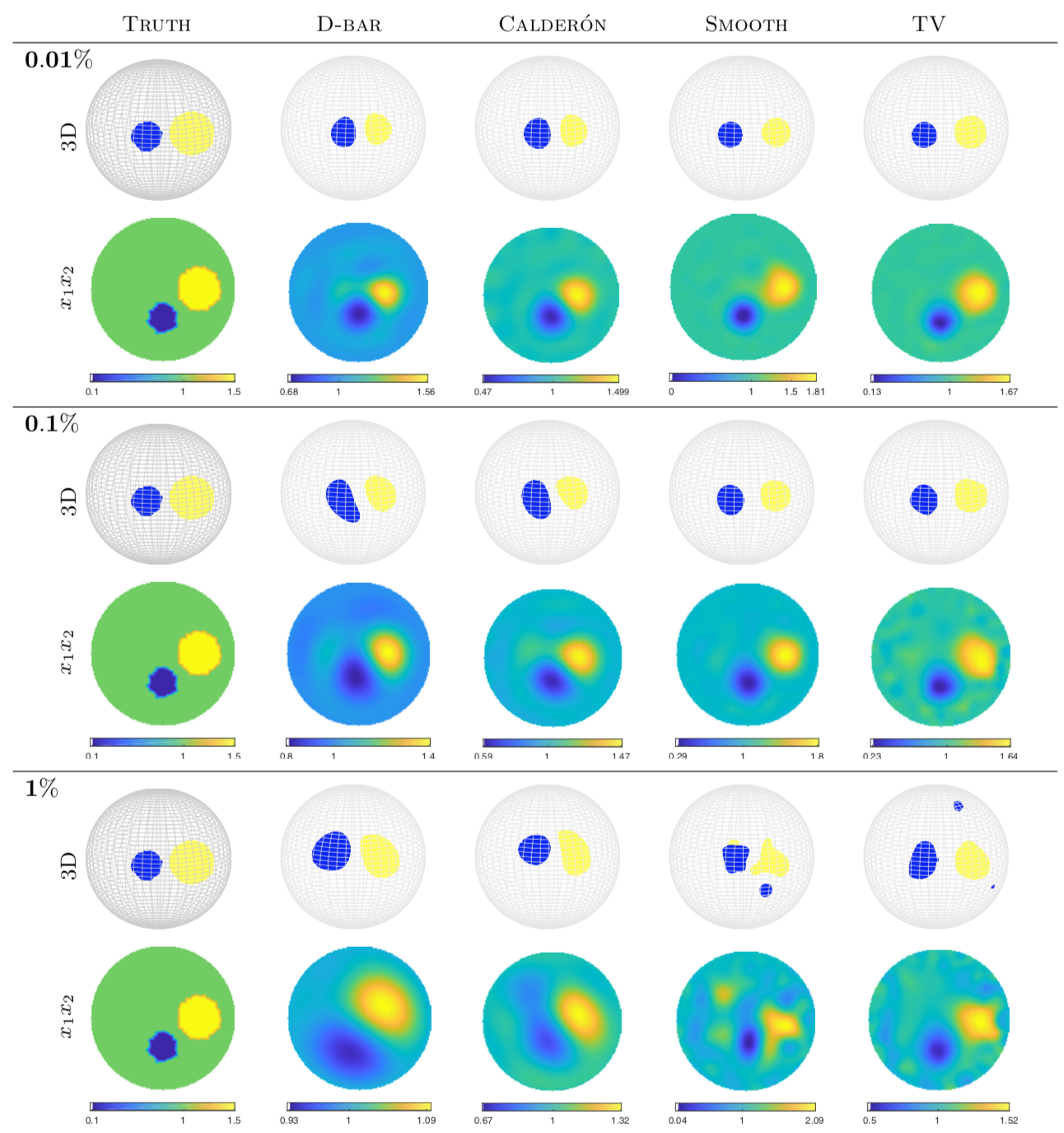}
\caption{\label{fig:T6_noisy} Comparison of reconstructions for the high contrast target T3 using various levels of noise.}
\end{figure}

How do the algorithms handle noisy voltage data?  We explore this question by adding $\eta*100\%$ Gaussian noise to the voltage data as described above in \S\ref{sec:methods-assess}.  We consider two target conductivities here, the real-valued heart and lungs target T2-B, and the high-contrast target T3.  Figures~\ref{fig:T3_noisy} and \ref{fig:T6_noisy} show the effect on targets T2-B and T3, respectively. The regularization parameters for Fig.~\ref{fig:T3_noisy} were: $\texp$ ($T_\xi =13,\;10, 7$), $\Fhat$ ($T_z=2.3,\;2, 1.3$, $t=0.05$), Smooth ($a = 0.444,\;0.888,\;0.888; \; b = 0.165,\;0.165,\;0.165; \; c = 10^{-8},\;10^{-8},\;10^{-8}$), TV ($\alpha_{\sigma} = 200,\; 400,\; 1000$, \; $\beta_{\sigma} = 10^{-8}$) for $0.01\%,\;0.1\%,$ and $1\%$ noise. Similarly, for Fig.~\ref{fig:T6_noisy} the parameters were: $\texp$ ($T_\xi =13,\;10, 7$), $\Fhat$ ($T_z=2,\;1.8, 1.4$, $t=0.05$), Smooth ($a = 0.444,\;0.888,\;0.888; \; b = 0.165,\;0.165,\;0.165; \; c = 10^{-8},\;10^{-8},\;10^{-8}$), TV ($\alpha_{\sigma} = 200,\; 400,\; 1000$, \; $\beta_{\sigma} = 10^{-8}$).  To show the 3D effect of noise more clearly, the top row at each noise level in Figures \ref{fig:T3_noisy} and \ref{fig:T6_noisy} shows the segmented targets used to compute the LE and RVR evaluation metrics.  Blue designates targets identified as resistors and yellow designates identified conductors.

\subsection{Evaluation Metrics}\label{ssec:eval_metrics_results}
In this section, we present the evaluation metric values, as defined in \S \ref{ssec:eval_metrics}, across all targets and methods. 
In the reporting of the metrics for T2, we label the lungs as ``lung 1" and ``lung 2," where ``lung 1" is the resistive target of larger volume. If the segmentation could not distinguish two separate lungs or identify clearly meaningful targets, the LE and RVR metrics relying on segmentation were not computed for those targets, and are denoted as N/A. If the segmentation identified ``targets" which were at least an order of magnitude smaller in volume than the true targets, those were omitted from the LE and RVR calculations.  As a reminder, dynamic range (DR), mean square error (MSE), and multi-scale structural similarity index (MS-SSIM) are computed for the whole image and not for individual targets.

Table \ref{tab:T7} shows the metrics for T1 reconstructions from simulated 128 electrode data corresponding to Figure \ref{fig:T7-compare}.   Table \ref{tab:T3cmplx} shows the metrics for T2-A reconstructions from simulated 128 electrode data corresponding to Figure \ref{fig:T3cmplx-compare}. The segmentation for the D-bar method identified two small ``targets," which were omitted from calculations. Likewise, one small ``target" was identified and omitted for the TV method. 

\clearpage


\begin{table}[h]
  \centering
  \footnotesize
  \caption{Evaluation metrics for T1 with 128 electrodes.}
    \begin{tabular}{|c||r|r|r|r|}
      \hline
                & \multicolumn{1}{c|}{D-bar} & \multicolumn{1}{c|}{Calder\'on} & \multicolumn{1}{c|}{Smooth} & \multicolumn{1}{c|}{TV} \\
    \hline
    \hline
    DR         & 95.25\% & 116.89\% & 162.62\% & 143.16\% \\
    \hline
    MSE        & 0.0305 & 0.0373 & 0.0141 & 0.0124 \\
    \hline
    MS-SSIM       & 0.8126 & 0.8324 & 0.8984 & 0.8978 \\
    \thickhline
    LE          & 0.0008 & 0.0008 & 0.0021 & 0.0007 \\
    \hline
    RVR         & 0.4435 & 0.3734 & 0.5088 & 0.6255 \\
    \hline
    \end{tabular}%
  \label{tab:T7}%
\end{table}%

\begin{table}[h]
  \centering
   \footnotesize
  \caption{Evaluation metrics for T2-A for 128 electrode data. Lung 1 is the resistive target with the largest volume.} 
    \begin{tabular}{|c|c|c||r|r|r|r|}
    \hline
       &       &       & \multicolumn{1}{c|}{D-bar} & \multicolumn{1}{c|}{Calder\'on} & \multicolumn{1}{c|}{Smooth} & \multicolumn{1}{c|}{TV} \\
    \hline
    \hline
    \multirow{2}[0]{*}{DR} & \multirow{2}[0]{*}{} & Re    & 50.19\% & 59.53\% & 92.59\% & 101.22\% \\
\cline{3-7}          &       & Im    & 98.02\% & 62.35\% & 85.80\% & 88.88\% \\
    \hline
    \multirow{2}[0]{*}{MSE} & \multirow{2}[0]{*}{} & Re    & 0.0097 & 0.0110 & 0.0042 & 0.0038 \\
\cline{3-7}          &       & Im    & 0.0022 & 0.0022 & 0.0011 & 0.0010 \\
    \hline
    \multirow{2}[0]{*}{MS-SSIM} & \multirow{2}[0]{*}{} & Re    & 0.8235 & 0.8193 & 0.7841 & 0.8179 \\
\cline{3-7}          &       & Im    & 0.8956 & 0.8710 & 0.8505 & 0.8558 \\
\thickhline
    \multirow{6}[0]{*}{LE} & \multirow{2}[0]{*}{ heart} & Re    & 0.1302 & 0.1172 & 0.0113 & 0.0263 \\
\cline{3-7}          &       & Im    & 0.1751 & 0.2133 & 0.1179 & 0.0949 \\
\cline{2-7}          & \multirow{2}[0]{*}{lung 1} & Re    & 0.1245 & 0.1401 & 0.0059 & 0.0105 \\
\cline{3-7}          &       & Im    & 0.1199 & 0.0843 & 0.0708 & 0.0412 \\
\cline{2-7}          & \multirow{2}[0]{*}{lung 2} & Re    & 0.1425 & 0.1746 & 0.0141 & 0.0067 \\
\cline{3-7}          &       & Im    & 0.1715 & 0.1290 & 0.0934 & 0.0538 \\
    \hline
    \multirow{6}[0]{*}{RVR} & \multirow{2}[0]{*}{ heart} & Re    & 0.8211 & 0.8112 & 1.0418 & 0.8497 \\
\cline{3-7}          &       & Im    & 0.2798 & 0.4819 & 1.0119 & 0.6845 \\
\cline{2-7}          & \multirow{2}[0]{*}{lung 1} & Re    & 0.9114 & 1.0072 & 1.4399 & 1.2159 \\
\cline{3-7}          &       & Im    & 1.1566 & 0.6418 & 0.9389 & 0.5171 \\
\cline{2-7}          & \multirow{2}[0]{*}{lung 2} & Re    & 0.8381 & 0.8041 & 1.4425 & 1.2289 \\
\cline{3-7}          &       & Im    & 1.3116 & 0.5547 & 0.5568 & 0.1976 \\
    \hline
    \end{tabular}%
  \label{tab:T3cmplx}%
\end{table}%

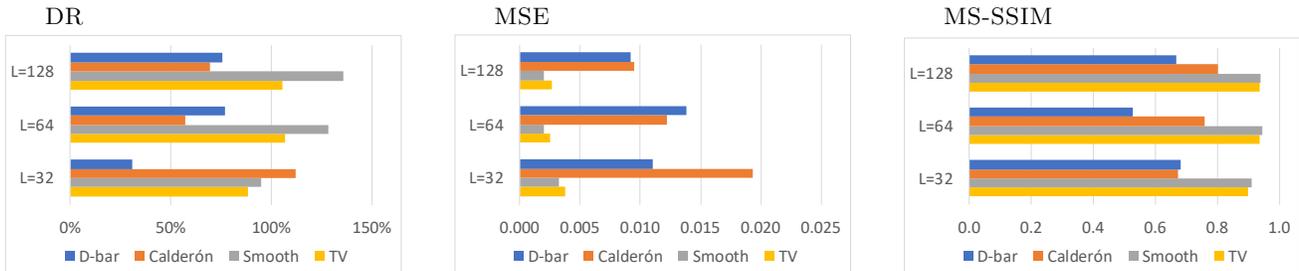
\begin{figure}[h!]
\centering
\begin{picture}(300,100)
\put(-100,0){\includegraphics[width=150pt]{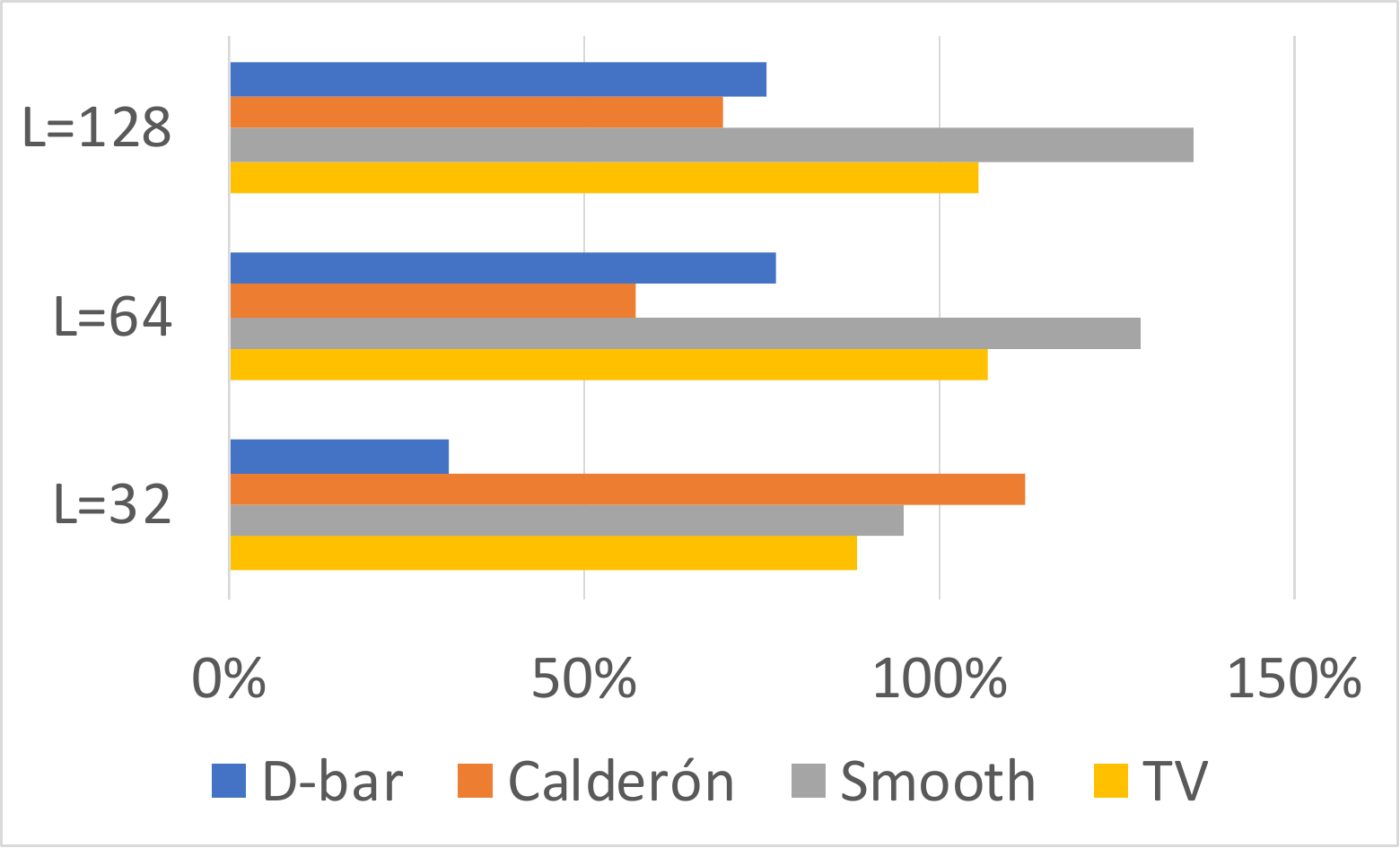}}
\put(70,0){\includegraphics[width=150pt]{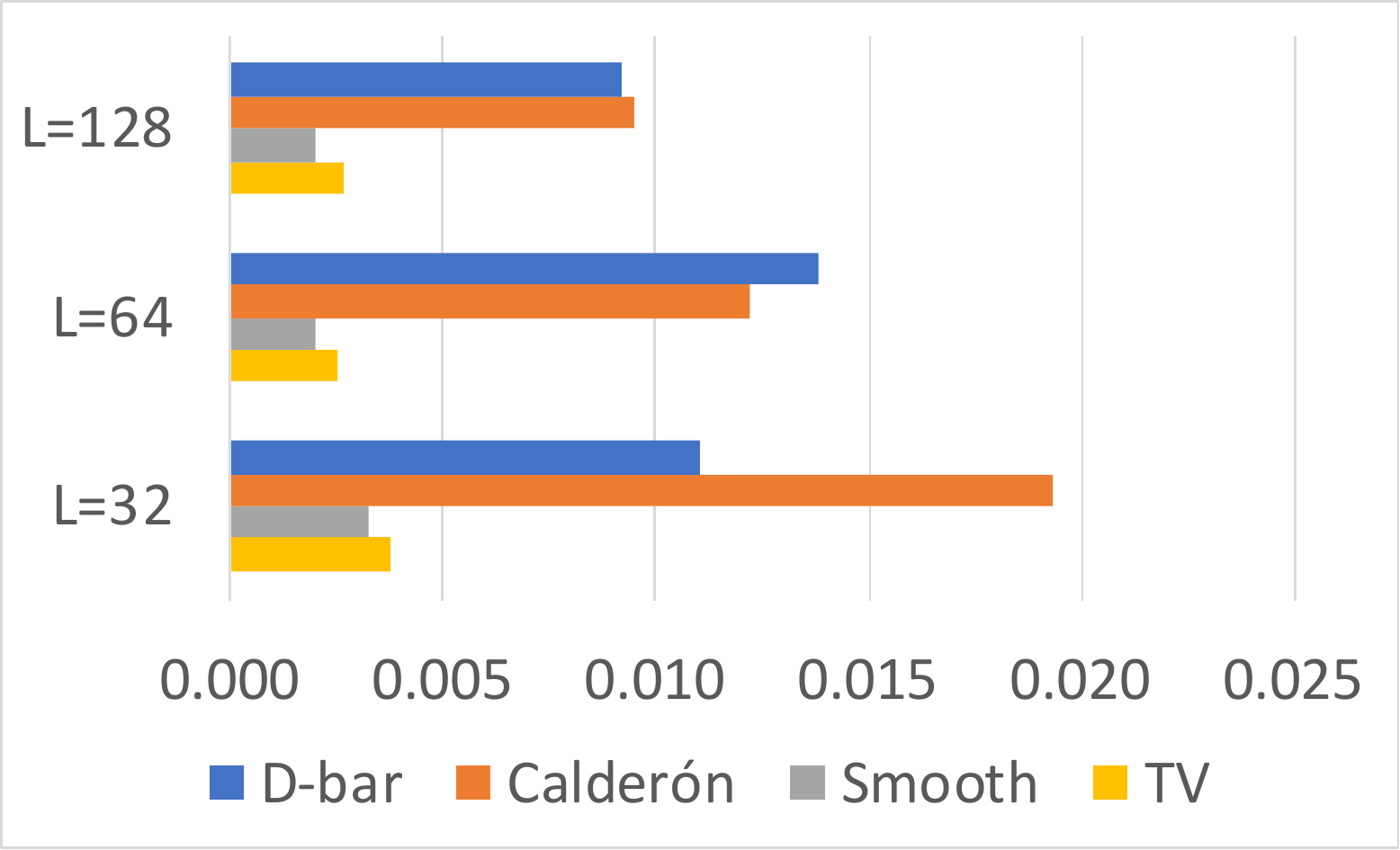}}
\put(240,0){\includegraphics[width=150pt]{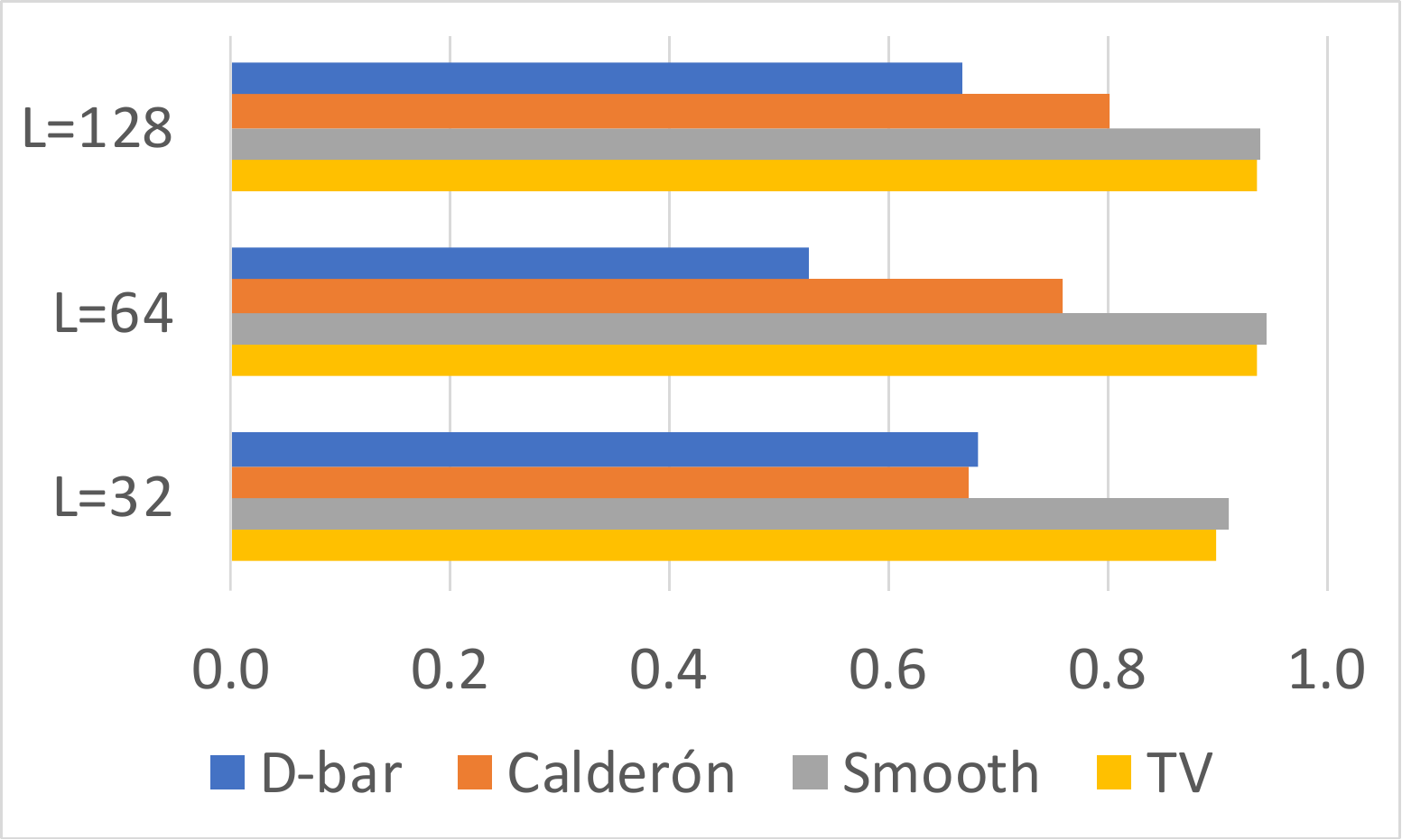}}

\put(-85,95){\footnotesize \sc DR}
\put(85,95){\footnotesize\sc MSE}
\put(255,95){\footnotesize \sc MS-SSIM}

\end{picture}
\caption{\label{fig:T3_elecs_metrics} Whole-image evaluation metrics for the real-valued heart and lungs target T2-B with decreasing numbers of simulated electrodes. Left: Dynamic Range, Middle: Mean Square Error, Right: Multi-Scale Structural Similarity Index.}
\end{figure}

Figure~\ref{fig:T3_elecs_metrics} and Table~\ref{tab:T3} show the metrics for T2-B reconstructions from simulated noise-free 128, 64 and 32 electrode data corresponding to Figure~\ref{fig:T3-compare}. Whole image metrics are compared in Figure~\ref{fig:T3_elecs_metrics}, while metrics relying on segmented targets are compared in Table ~\ref{tab:T3}. 
\vspace{0.5em}
\begin{table}[h]
  \centering
  \footnotesize
  \caption{T2-B evaluation metrics across all electrode configurations considered. Lung 1 is the resistive target with the larger volume.}
    \begin{tabular}{|c|c|c||c|c|r|r|}
    \hline
    \multicolumn{1}{|r|}{} &       &       & D-bar & Calder\'on & \multicolumn{1}{c|}{Smooth} & \multicolumn{1}{c|}{TV} \\
    \hline
    \hline
    \multirow{9}[0]{*}{LE} & \multirow{3}[0]{*}{heart} & L=128 & \multicolumn{1}{r|}{0.1545} & \multicolumn{1}{r|}{0.1296} & 0.0086 & 0.0171 \\
\cline{3-7}          &       & L=64  & \multicolumn{1}{r|}{0.2416} & \multicolumn{1}{r|}{0.1899} & 0.0070 & 0.0176 \\
\cline{3-7}          &       & L=32  & \multicolumn{1}{r|}{0.2424} & \multicolumn{1}{r|}{0.1837} & 0.0115 & 0.0035 \\
\clineB{2-7}{2.5}          & \multirow{3}[0]{*}{lung 1} & L=128 & \multicolumn{1}{r|}{0.1099} & \multicolumn{1}{r|}{0.1216} & 0.0130 & 0.0194 \\
\cline{3-7}          &       & L=64  & \multicolumn{1}{r|}{0.1627} & \multicolumn{1}{r|}{0.1628} & 0.0084 & 0.0118 \\
\cline{3-7}          &       & L=32  & N/A   & N/A   & 0.0040 & 0.0098 \\
\clineB{2-7}{2.5}          & \multirow{3}[0]{*}{lung 2} & L=128 & \multicolumn{1}{r|}{0.1171} & \multicolumn{1}{r|}{0.1531} & 0.0144 & 0.0109 \\
\cline{3-7}          &       & L=64  & \multicolumn{1}{r|}{0.2221} & \multicolumn{1}{r|}{0.2172} & 0.0069 & 0.0085 \\
\cline{3-7}          &       & L=32  & N/A   & N/A   & 0.0135 & 0.0068 \\
    \hline
    \hline
    \multirow{9}[0]{*}{RVR} & \multirow{3}[0]{*}{heart} & L=128 & \multicolumn{1}{r|}{0.7526} & \multicolumn{1}{r|}{1.1204} & 0.6520 & 0.9152 \\
\cline{3-7}         &       & L=64  & \multicolumn{1}{r|}{1.0722} & \multicolumn{1}{r|}{1.1101} & 0.7737 & 0.9869 \\
\cline{3-7}          &       & L=32  & \multicolumn{1}{r|}{2.9867} & \multicolumn{1}{r|}{1.9449} & 1.2085 & 1.2273 \\
\clineB{2-7}{2.5}          & \multirow{3}[0]{*}{lung 1} & L=128 & \multicolumn{1}{r|}{0.3977} & \multicolumn{1}{r|}{0.6102} & 0.6360 & 0.7102 \\
\cline{3-7}          &       & L=64  & \multicolumn{1}{r|}{0.3389} & \multicolumn{1}{r|}{0.5475} & 0.7038 & 0.8359 \\
\cline{3-7}          &       & L=32  & N/A   & N/A   & 0.9249 & 0.9318 \\
\clineB{2-7}{2.5}          & \multirow{3}[0]{*}{lung 2} & L=128 & \multicolumn{1}{r|}{0.3312} & \multicolumn{1}{r|}{0.2990} & 0.6342 & 0.6600 \\
\cline{3-7}          &       & L=64  & \multicolumn{1}{r|}{0.2311} & \multicolumn{1}{r|}{0.3047} & 0.6908 & 0.8039 \\
\cline{3-7}          &       & L=32  & N/A   & N/A   & 0.8597 & 0.8590 \\
    \hline
    \end{tabular}%
  \label{tab:T3}%
\end{table}%

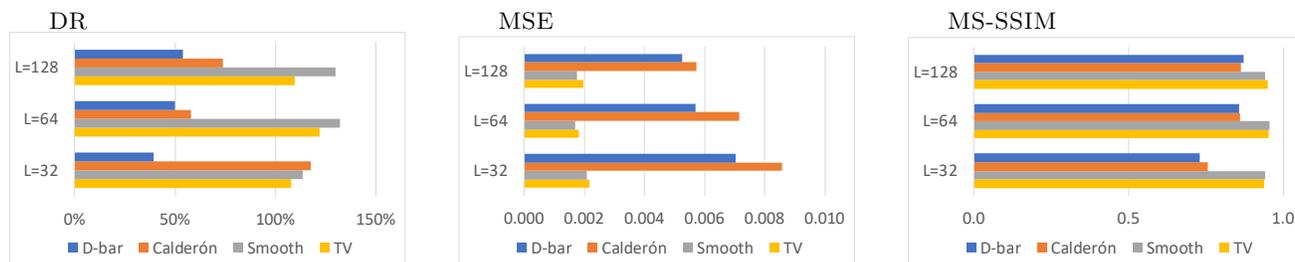
\begin{figure}[h!]
\centering
\begin{picture}(300,100)
\put(-100,0){\includegraphics[width=150pt]{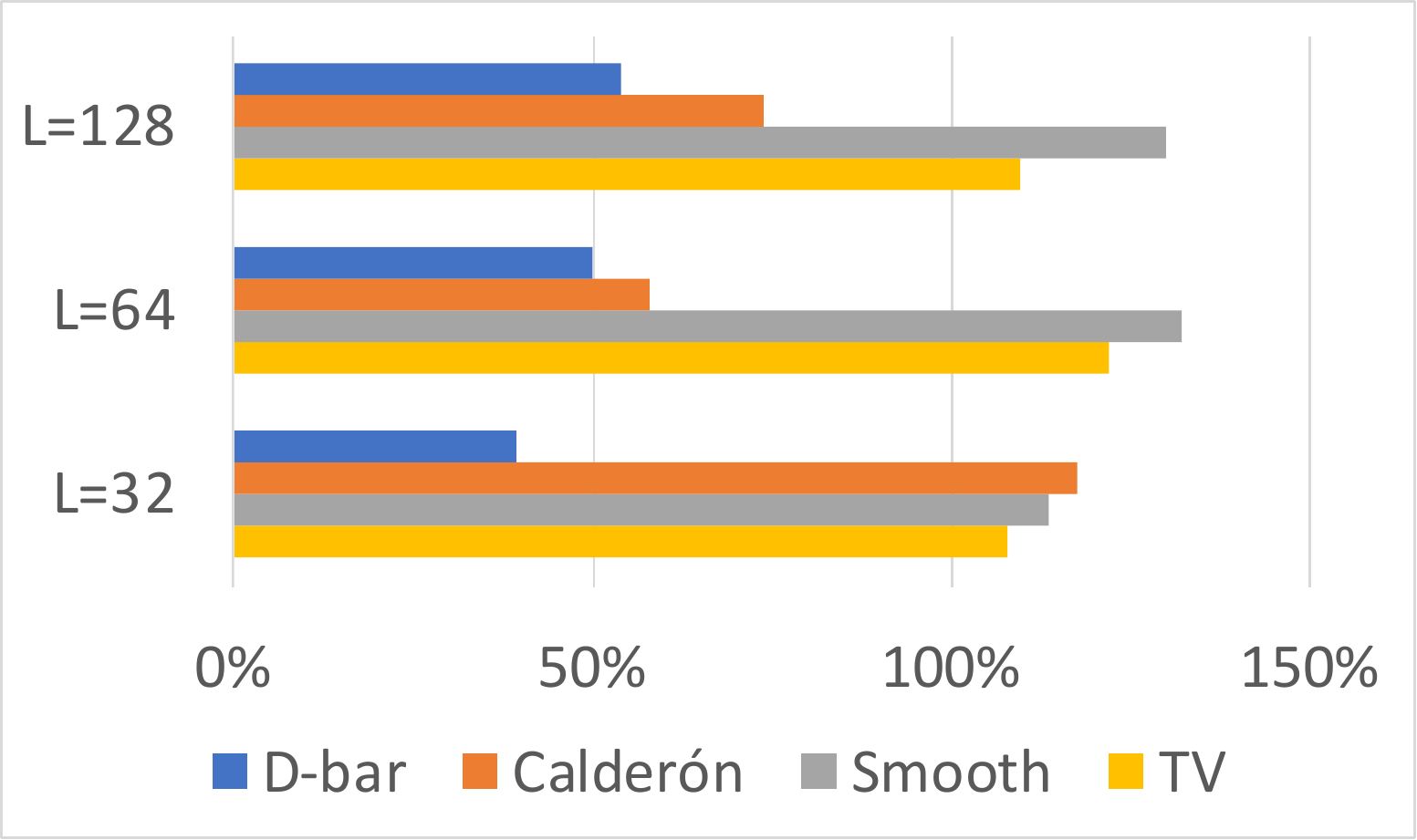}}
\put(70,0){\includegraphics[width=150pt]{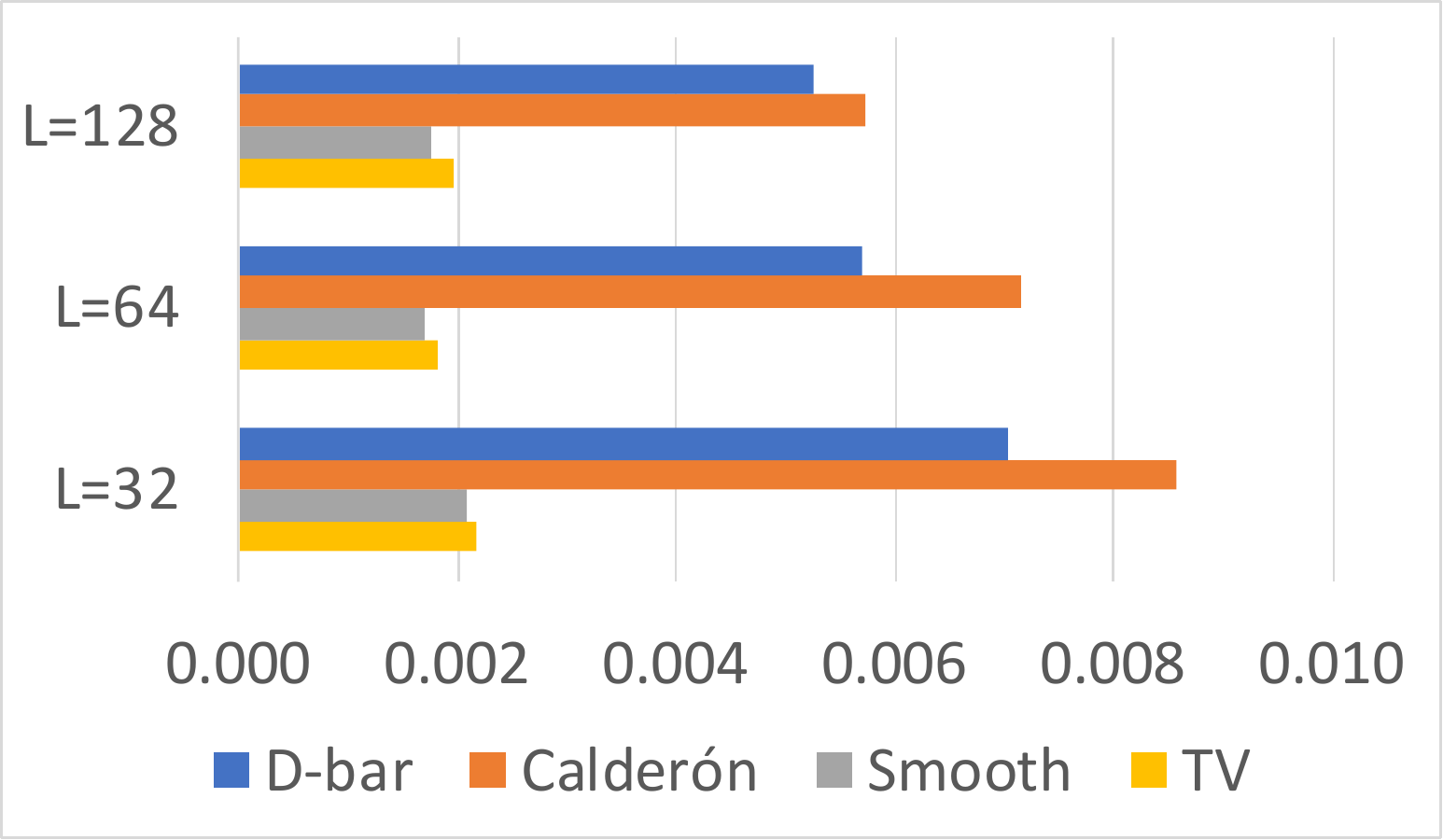}}
\put(240,0){\includegraphics[width=150pt]{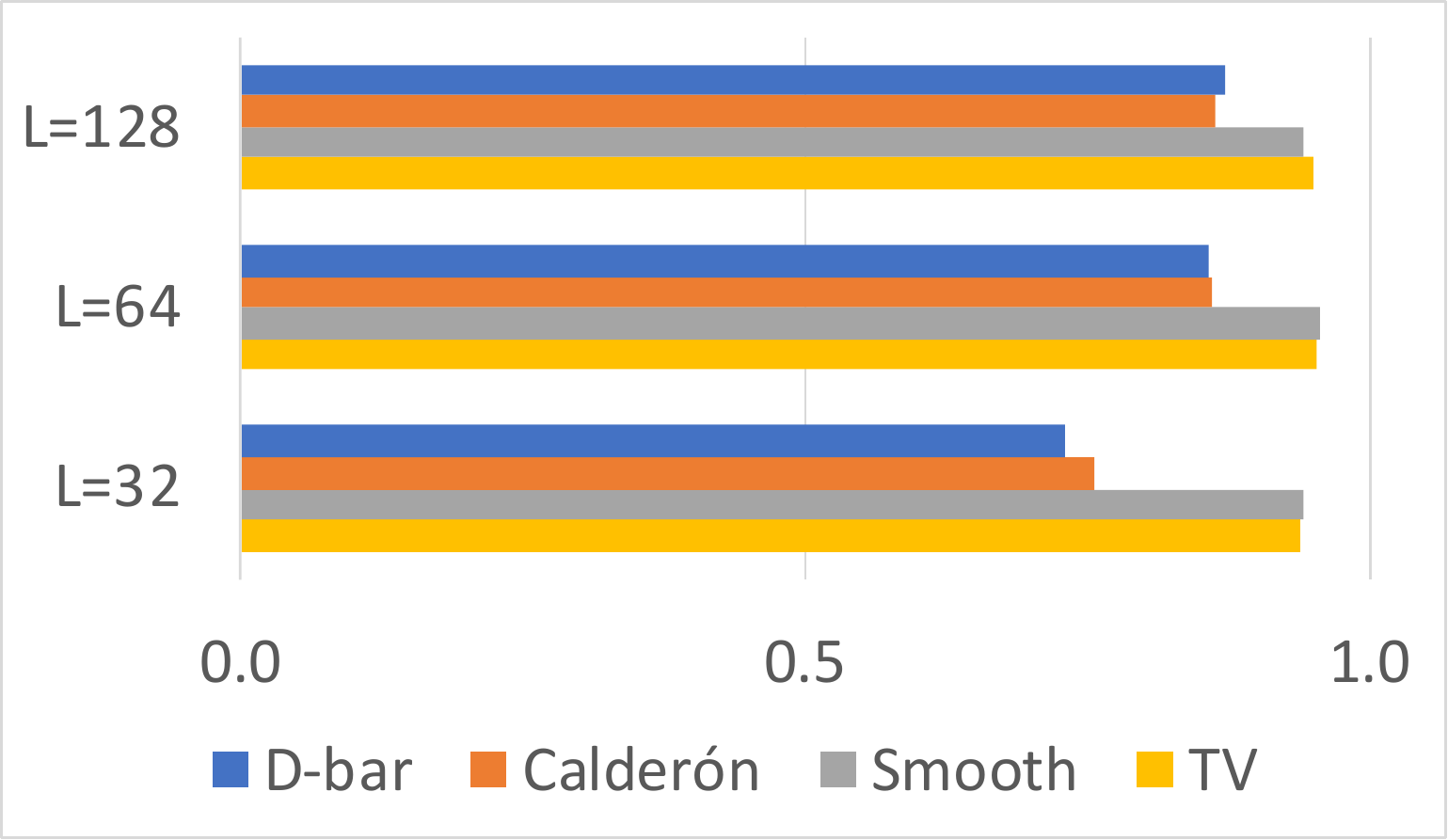}}

\put(-85,90){\footnotesize \sc DR}
\put(85,90){\footnotesize\sc MSE}
\put(255,90){\footnotesize \sc MS-SSIM}

\end{picture}
\caption{\label{fig:T6_elecs_metrics} Whole-image evaluation metrics for target T3 with decreasing numbers of simulated electrodes. Left: Dynamic Range, Middle: Mean Square Error, Right: Multi-Scale Structural Similarity Index}
\end{figure}

Figure~\ref{fig:T6_elecs_metrics} and Table~\ref{tab:T6} show the metrics for T2-B reconstructions from simulated noise-free 128, 64 and 32 electrode data corresponding to Figure~\ref{fig:T6-compare}. Whole image metrics are compared in Figure~\ref{fig:T6_elecs_metrics}, while metrics relying on segmented targets are compared in Table ~\ref{tab:T6}.
\vspace{0.5em}
\begin{table}[h!]
  \centering
  \footnotesize
  \caption{Evaluation metrics for the high-contrast example T3 across all electrode configurations considered.}
    \begin{tabular}{|c|c|c||r|r|r|r|}
    \hline
    \multicolumn{1}{|r|}{} &       &       & \multicolumn{1}{c|}{D-bar} & \multicolumn{1}{c|}{Calderon} & \multicolumn{1}{c|}{Smooth} & \multicolumn{1}{c|}{TV} \\
    \hline
    \hline
    \multirow{6}[0]{*}{LE} & \multirow{3}[0]{*}{ conductor} & L=128 & 0.1267 & 0.1181 & 0.0264 & 0.0191 \\
\cline{3-7}          &       & L=64  & 0.1522 & 0.1566 & 0.0085 & 0.0136 \\
\cline{3-7}          &       & L=32  & 0.1559 & 0.1548 & 0.0111 & 0.0093 \\
\clineB{2-7}{2.5}          & \multirow{3}[6]{*}{ resistor} & L=128 & 0.0853 & 0.0861 & 0.0128 & 0.0074 \\
\cline{3-7}          &       & L=64  & 0.1072 & 0.1107 & 0.0050 & 0.0061 \\
\cline{3-7}          &       & L=32  & 0.1465 & 0.1062 & 0.0057 & 0.0109 \\
    \hline
    \hline
    \multirow{6}[0]{*}{RVR} & \multirow{3}[0]{*}{ conductor} & L=128 & 0.4672 & 0.4690 & 0.4122 & 0.6124 \\
\cline{3-7}          &       & L=64  & 0.5357 & 0.4343 & 0.4456 & 0.5490 \\
\cline{3-7}          &       & L=32  & 0.5883 & 0.5580 & 0.5492 & 0.6912 \\
\clineB{2-7}{2.5}          & \multirow{3}[0]{*}{ resistor} & L=128 & 1.6747 & 1.6273 & 0.9951 & 1.0290 \\
\cline{3-7}          &       & L=64  & 1.9706 & 1.6022 & 1.0895 & 1.0479 \\
\cline{3-7}          &       & L=32  & 2.1047 & 2.0931 & 1.4676 & 1.1714 \\
    \hline
    \end{tabular}%
  \label{tab:T6}%
\end{table}%

Figure~\ref{fig:T3_noisy_metrics} and Table~\ref{tab:T3_noisy} show the metrics for T2-B reconstructions from simulated 128 electrode data with 0.01\%, 0.1\%, and 1\% noise corresponding to Figure~\ref{fig:T3_noisy}. Whole image metrics are compared in Figure~\ref{fig:T3_noisy_metrics}, while metrics relying on segmented targets are compared in Table ~\ref{tab:T3_noisy}. For the TV method at 1\% noise, the segmentation identified 12 small ``targets", which were omitted from calculations. At 1\% noise, the conductive threshold for the D-bar method had to be set to 0.52. At 1\% noise segmentation failed to produce identifiable targets for the Smooth LS method.
\begin{figure}[h!]
\centering
\begin{picture}(300,100)
\put(-100,0){\includegraphics[width=150pt]{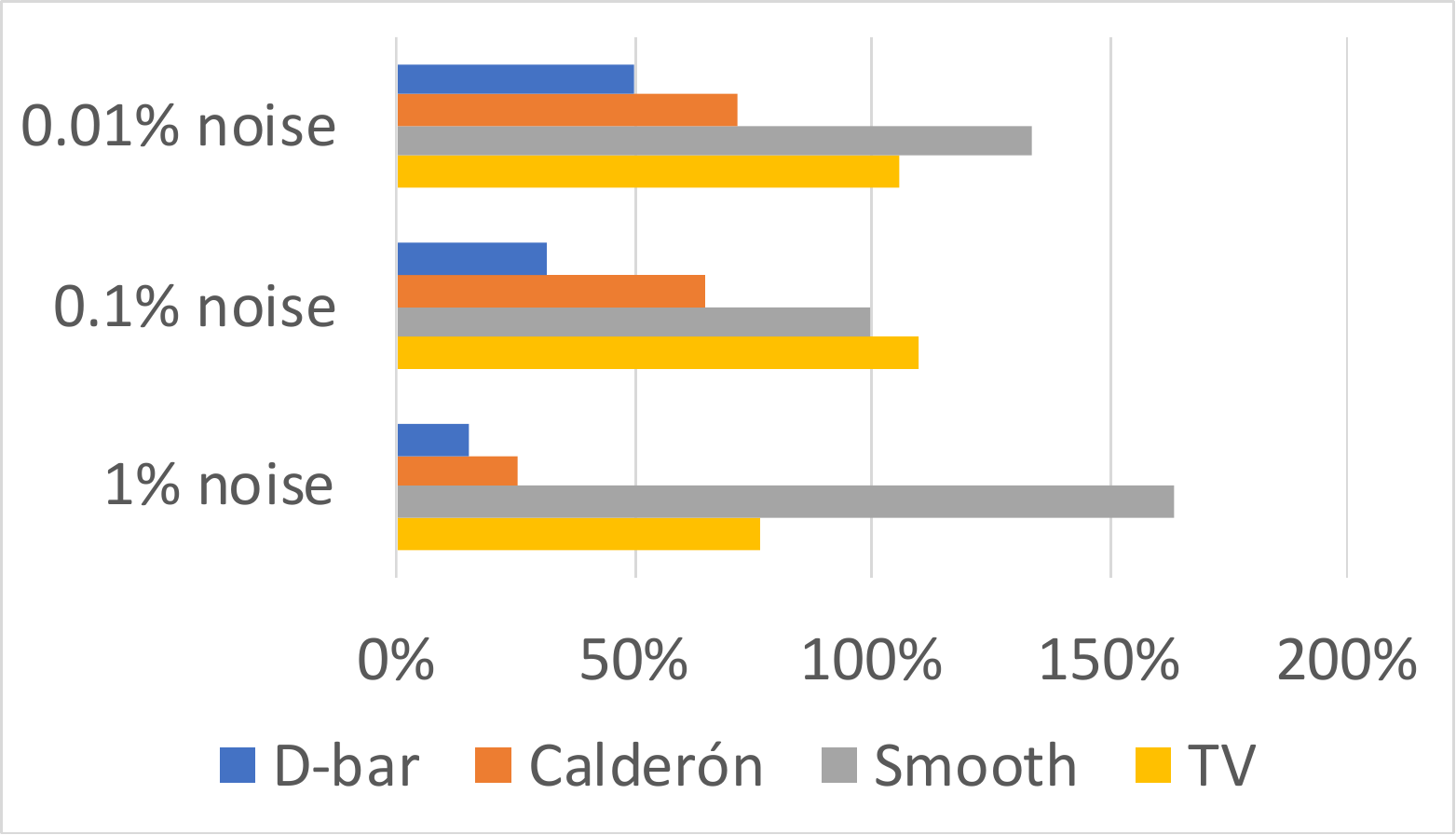}}
\put(70,0){\includegraphics[width=150pt]{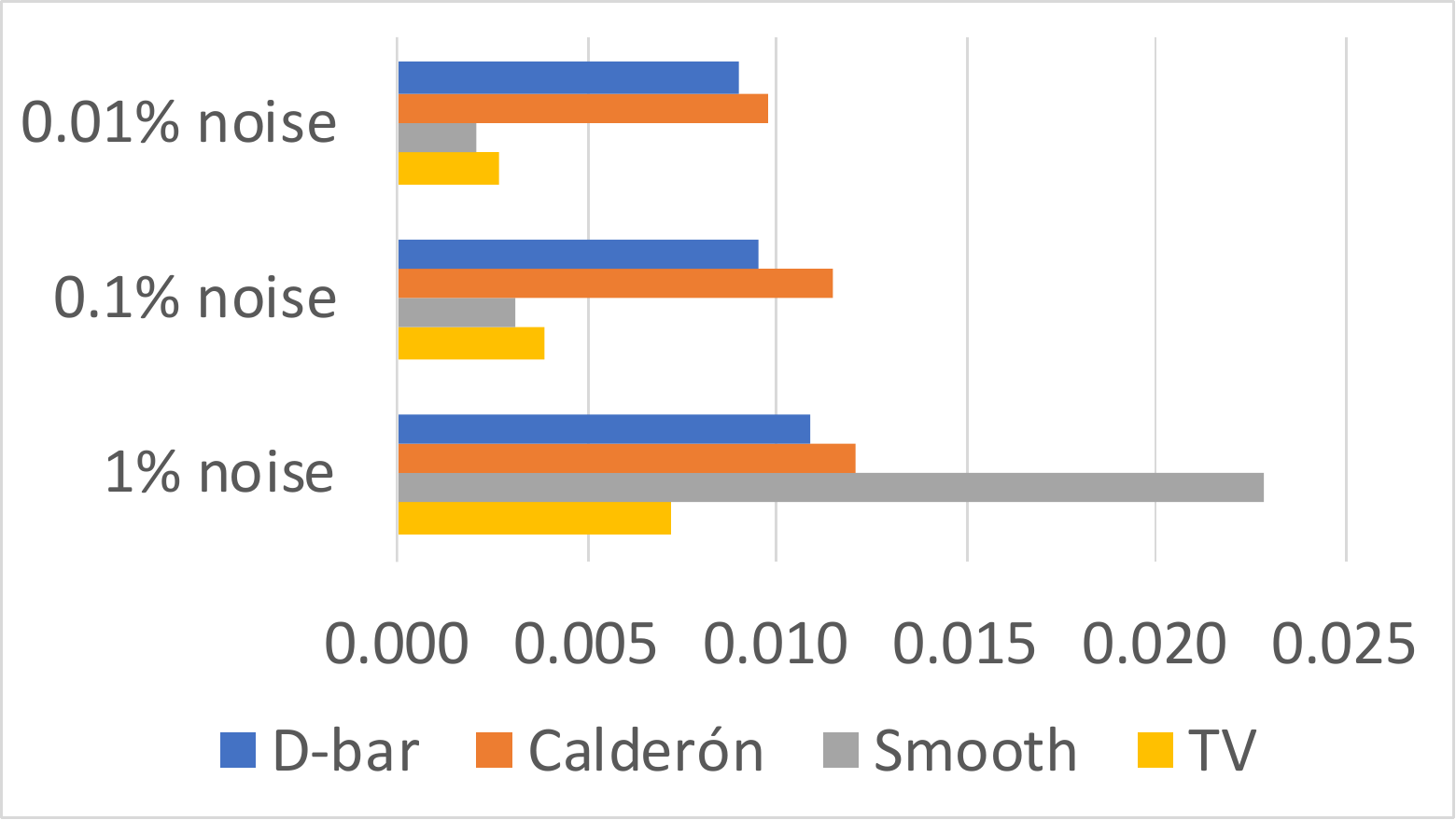}}
\put(240,0){\includegraphics[width=150pt]{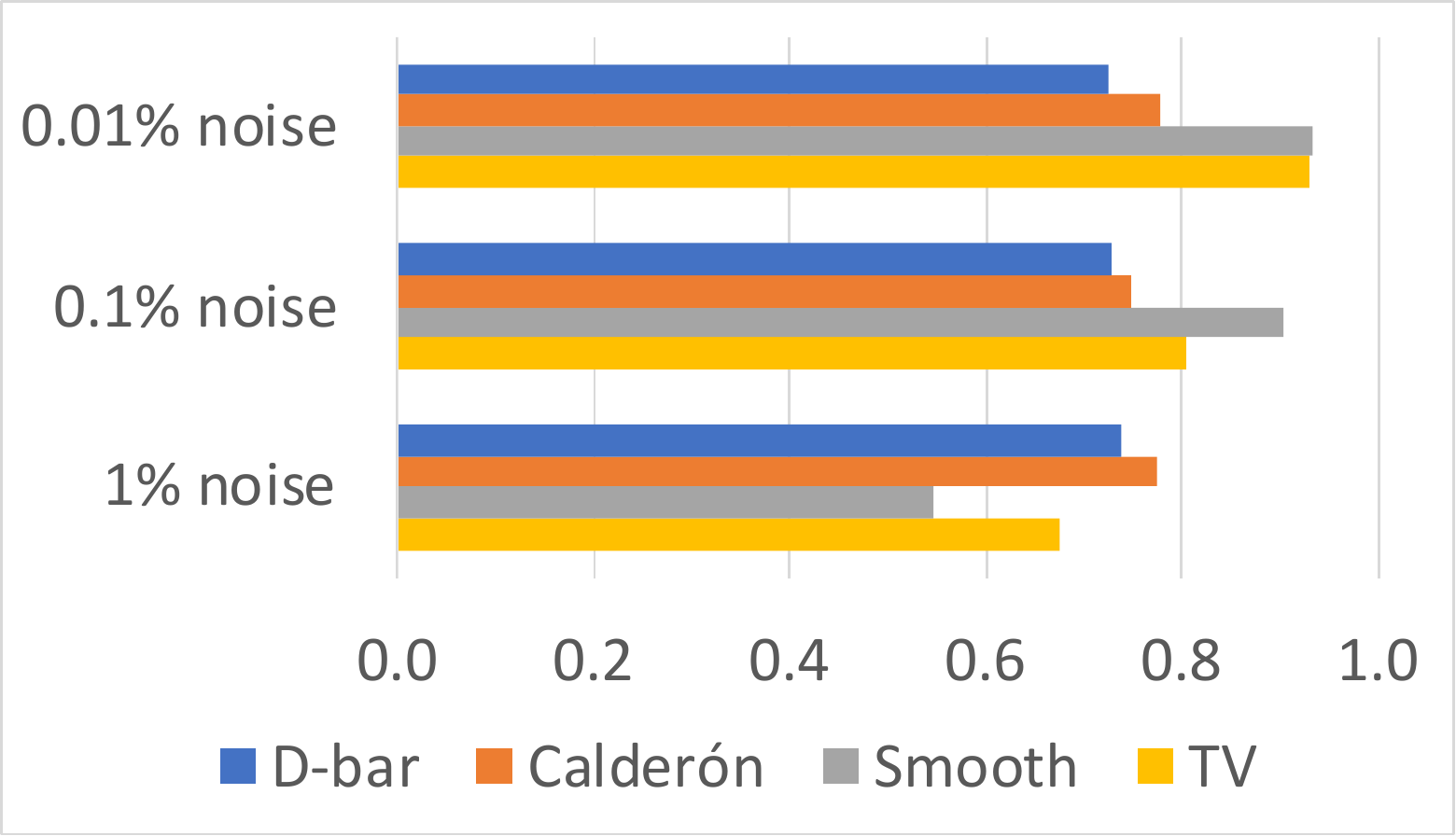}}

\put(-85,90){\footnotesize \sc DR}
\put(85,90){\footnotesize\sc MSE}
\put(255,90){\footnotesize \sc MS-SSIM}

\end{picture}
\caption{\label{fig:T3_noisy_metrics} Whole-image evaluation metrics for the real-valued heart and lungs target T2-B with increasing levels of noise added to the voltage data. Left: Dynamic Range, Middle: Mean Square Error, Right: Multi-Scale Structural Similarity Index}
\end{figure}
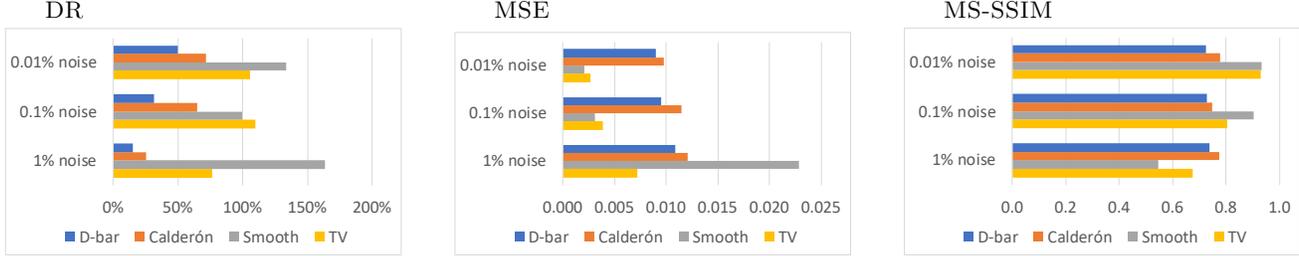
\begin{table}[h]
  \centering
  \footnotesize
  \caption{Evaluation metrics for T2-B with 0.01\%, 0.1\%, and 1\% noise with $128$ electrodes.}
    \begin{tabular}{|c|c|c||c|c|c|r|}
    \hline
          &       &       & D-bar & Calderon & Smooth & \multicolumn{1}{c|}{TV} \\
    \hline
    \hline
    \multirow{9}[0]{*}{LE} & \multirow{3}[0]{*}{ heart} & 0.01\% noise & \multicolumn{1}{r|}{0.1451} & \multicolumn{1}{r|}{0.1494} & \multicolumn{1}{r|}{0.0096} & 0.0158 \\
\cline{3-7}          &       & 0.1\% noise & \multicolumn{1}{r|}{0.2372} & \multicolumn{1}{r|}{0.2147} & \multicolumn{1}{r|}{0.0184} & 0.0286 \\
\cline{3-7}          &       & 1\% noise & \multicolumn{1}{r|}{0.2061} & \multicolumn{1}{r|}{0.1520} & N/A   & 0.0770 \\
\clineB{2-7}{2.5}          & \multirow{3}[0]{*}{lung 1} & 0.01\% noise & \multicolumn{1}{r|}{0.1294} & \multicolumn{1}{r|}{0.1201} & \multicolumn{1}{r|}{0.0138} & 0.0209 \\
\cline{3-7}          &       & 0.1\% noise & \multicolumn{1}{r|}{0.1541} & \multicolumn{1}{r|}{0.1540} & \multicolumn{1}{r|}{0.0125} & 0.0213 \\
\cline{3-7}          &       & 1\% noise & N/A   & N/A   & N/A   & 0.0513 \\
\clineB{2-7}{2.5}          & \multirow{3}[0]{*}{lung 2} & 0.01\% noise & \multicolumn{1}{r|}{0.1097} & \multicolumn{1}{r|}{0.1336} & \multicolumn{1}{r|}{0.0141} & 0.0111 \\
\cline{3-7}          &       & 0.1\% noise & \multicolumn{1}{r|}{0.1245} & \multicolumn{1}{r|}{0.1498} & \multicolumn{1}{r|}{0.0210} & 0.0220 \\
\cline{3-7}          &       & 1\% noise & N/A   & N/A   & N/A   & 0.0444 \\
    \hline
    \hline
    \multirow{9}[0]{*}{RVR} & \multirow{3}[0]{*}{ heart} & 0.01\% noise & \multicolumn{1}{r|}{1.3702} & \multicolumn{1}{r|}{1.0072} & \multicolumn{1}{r|}{0.6605} & 0.9279 \\
\cline{3-7}          &       & 0.1\% noise & \multicolumn{1}{r|}{2.5372} & \multicolumn{1}{r|}{1.3738} & \multicolumn{1}{r|}{1.0518} & 1.1021 \\
\cline{3-7}          &       & 1\% noise & \multicolumn{1}{r|}{0.6135} & \multicolumn{1}{r|}{3.7011} & \multicolumn{1}{l|}{N/A} & 1.0894 \\
\clineB{2-7}{2.5}          & \multirow{3}[0]{*}{lung 1} & 0.01\% noise & \multicolumn{1}{r|}{0.5548} & \multicolumn{1}{r|}{0.6647} & \multicolumn{1}{r|}{0.6603} & 0.7192 \\
\cline{3-7}          &       & 0.1\% noise & \multicolumn{1}{r|}{1.0583} & \multicolumn{1}{r|}{0.7483} & \multicolumn{1}{r|}{0.9256} & 0.4780 \\
\cline{3-7}          &       & 1\% noise & N/A   & N/A   & N/A   & 0.9519 \\
\clineB{2-7}{2.5}          & \multirow{3}[0]{*}{lung 2} & 0.01\% noise & \multicolumn{1}{r|}{0.4064} & \multicolumn{1}{r|}{0.3625} & \multicolumn{1}{r|}{0.6504} & 0.6581 \\
\cline{3-7}          &       & 0.1\% noise & \multicolumn{1}{r|}{1.1779} & \multicolumn{1}{r|}{0.5656} & \multicolumn{1}{r|}{0.9098} & 0.4665 \\
\cline{3-7}          &       & 1\% noise & N/A   & N/A   & N/A   & 0.7636 \\
\hline    \end{tabular}%
  \label{tab:T3_noisy}%
\end{table}%


\vspace{0.5em}

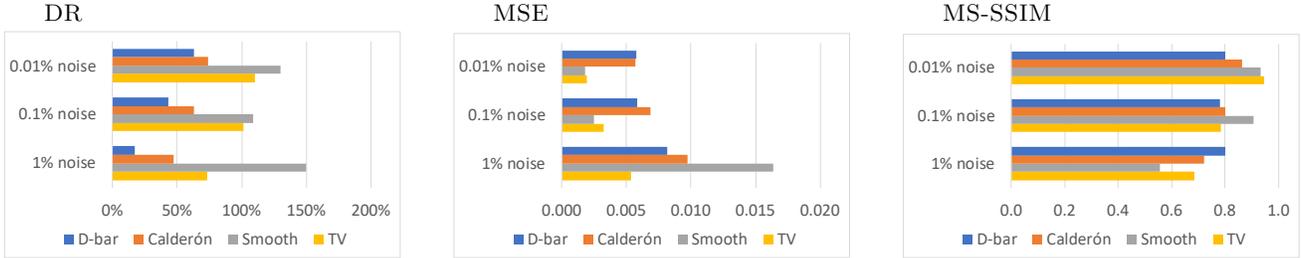
\begin{figure}[h!]
\centering
\begin{picture}(300,100)
\put(-100,0){\includegraphics[width=150pt]{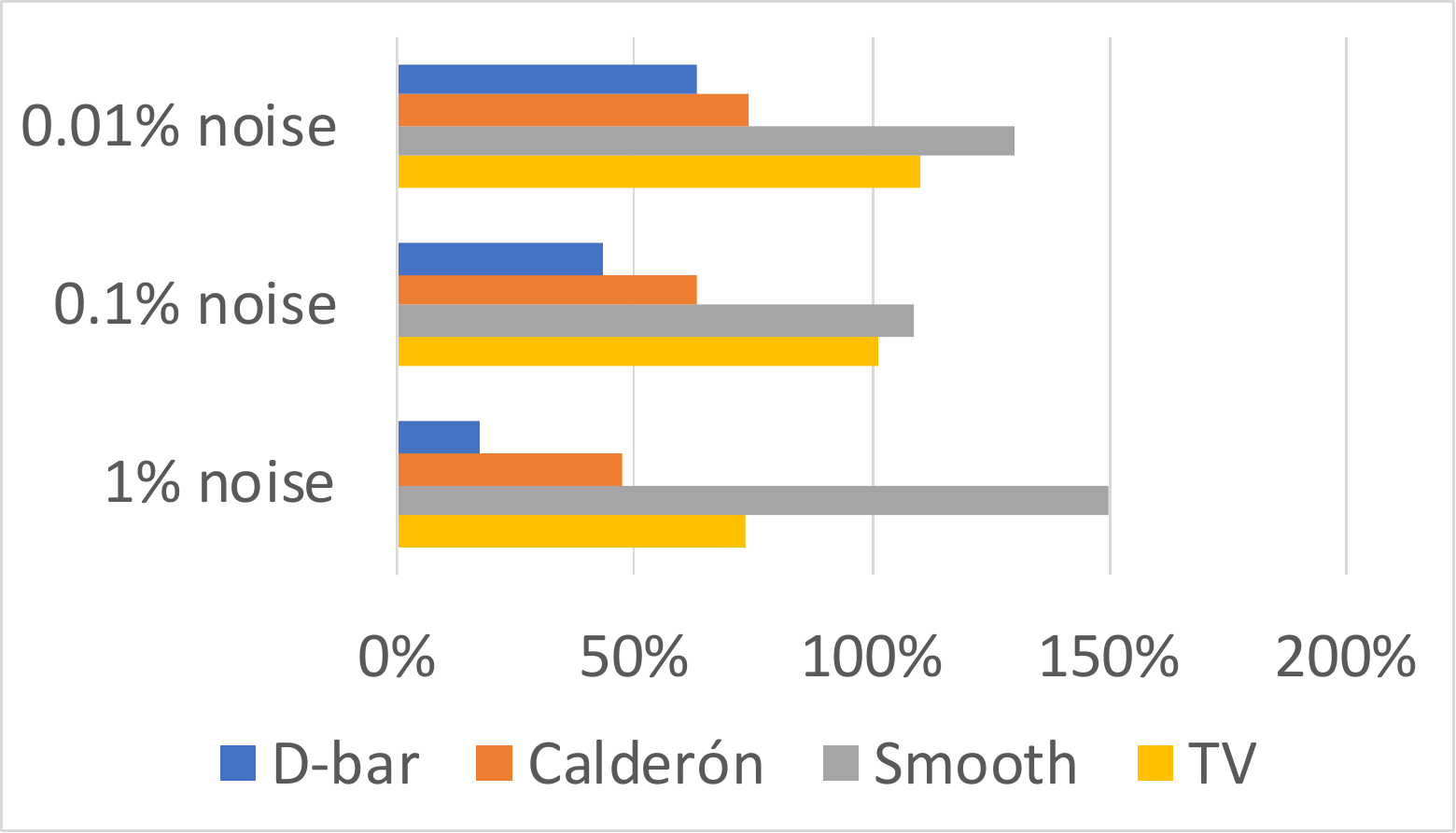}}
\put(70,0){\includegraphics[width=150pt]{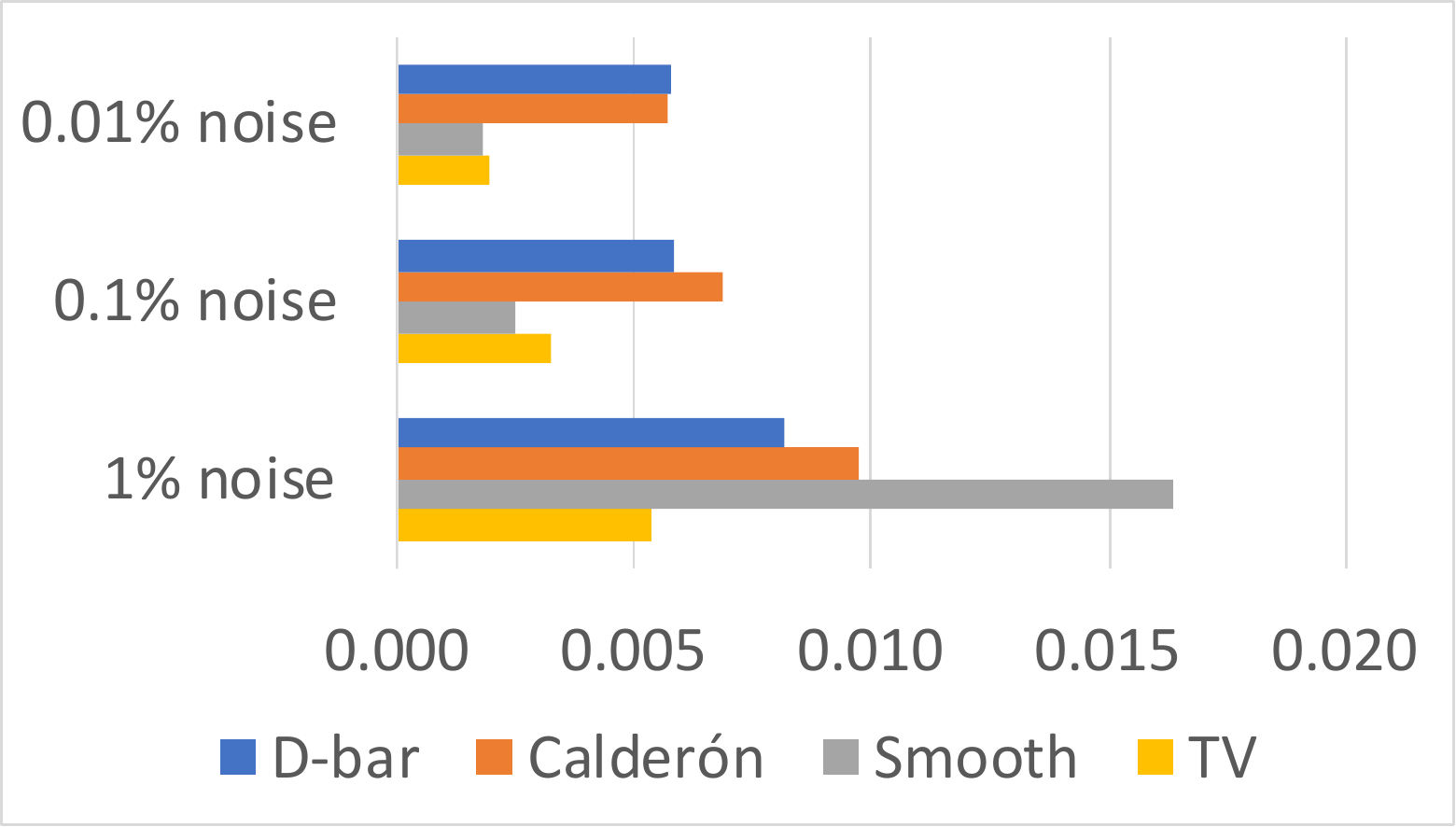}}
\put(240,0){\includegraphics[width=150pt]{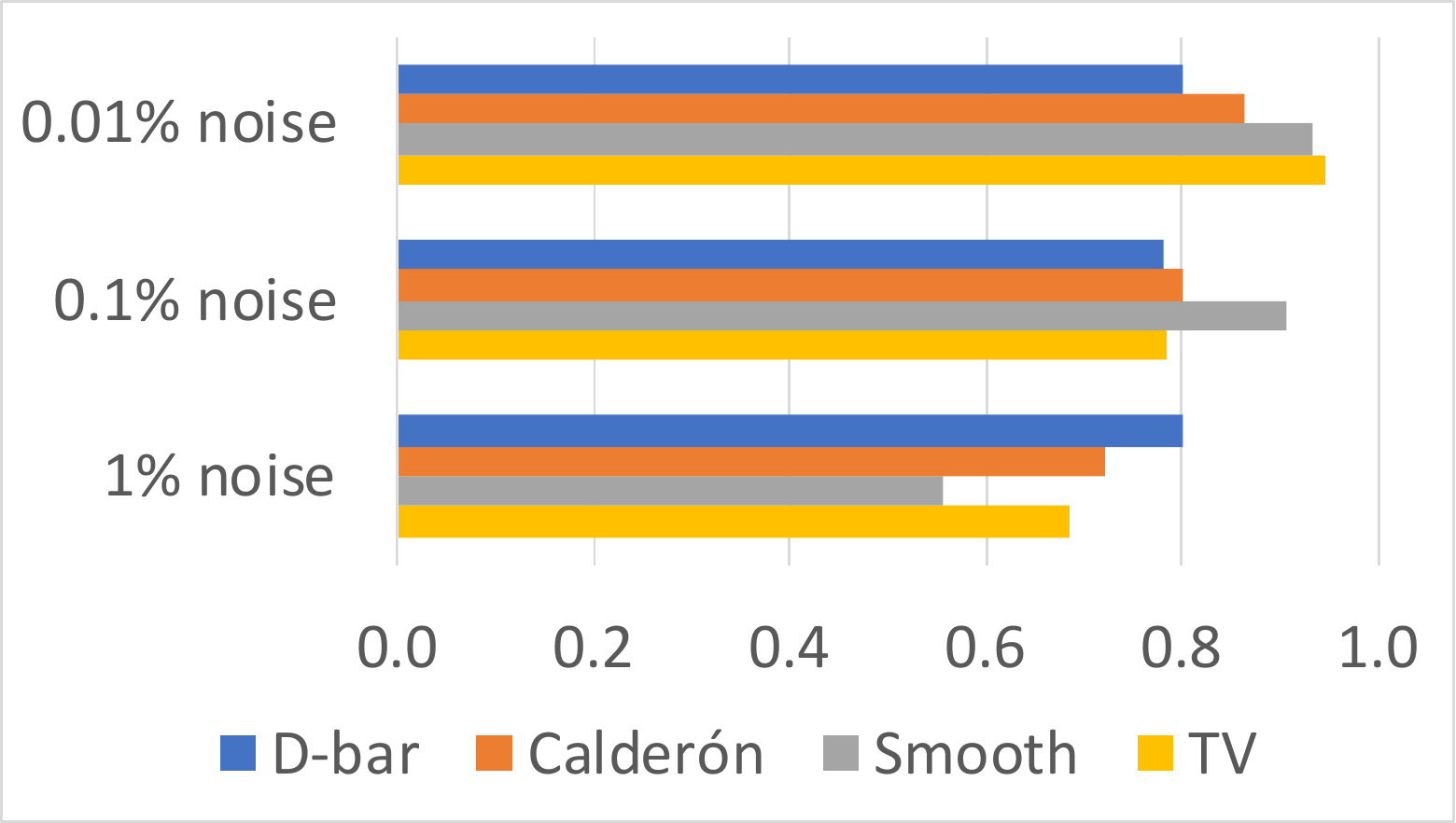}}

\put(-85,90){\footnotesize \sc DR}
\put(85,90){\footnotesize\sc MSE}
\put(255,90){\footnotesize \sc MS-SSIM}

\end{picture}
\caption{\label{fig:T6_noisy_metrics} Whole-image evaluation metrics for the target T3 with increasing levels of noise added to the voltage data. Left: Dynamic Range, Middle: Mean Square Error, Right: Multi-Scale Structural Similarity Index}
\end{figure}

Figure~\ref{fig:T6_noisy_metrics} and Table~\ref{tab:T6_noisy} show the metrics for the high contrast T3 reconstructions from simulated 128 electrode data with 0.01\%, 0.1\%, and 1\% noise corresponding to Figure~\ref{fig:T6_noisy}. Whole image metrics are compared in Figure~\ref{fig:T6_noisy_metrics}, while metrics relying on segmented targets are compared in Table ~\ref{tab:T6_noisy}. We note that in the 1\% noise the conductive and resistive thresholds were changed from 0.5 for the CGO methods. For D-bar, the conductive and resistive thresholds were set to 0.4 and 0.6, respectively. For Calder\'on, both were set to 0.6. For the TV method at 1\% noise, the segmentation identified five small ``targets," which were omitted from calculations. For the Smooth LS method, the segmentation found three small ``targets" at 1\% noise, which were omitted from calculations\footnote{One of these ``targets" was only 8 times smaller in volume than the resistive target, but is clearly not an intended reconstructed target by inspection.}.

\begin{table}[h]
  \centering
  \footnotesize
  \caption{Evaluation metrics for the high-contrast example T3 with 0.01\%, 0.1\%, and 1\% noise with $128$ electrodes.}
    \begin{tabular}{|c|c|c||r|r|r|r|}
    \hline
          &       &       & \multicolumn{1}{c|}{D-bar} & \multicolumn{1}{c|}{Calderon} & \multicolumn{1}{c|}{Smooth} & \multicolumn{1}{c|}{TV} \\
    \hline
    \hline
    \multirow{6}[0]{*}{LE} & \multirow{3}[0]{*}{ conductor} & 0.01\% noise & 0.1303 & 0.1155 & 0.0292 & 0.0204 \\
\cline{3-7} &   & 0.1\% noise & 0.1044 & 0.1391 & 0.0127 & 0.0262 \\
\cline{3-7}   &       & 1\% noise & 0.2501 & 0.1854 & 0.0611 & 0.0551 \\
\clineB{2-7}{2.5}     & \multirow{3}[0]{*}{ resistor} & 0.01\% noise & 0.0887 & 0.0814 & 0.0134 & 0.0080 \\
\cline{3-7}    &       & 0.1\% noise & 0.1444 & 0.1046 & 0.0171 & 0.0046 \\
\cline{3-7}     &       & 1\% noise & 0.2193 & 0.1930 & 0.1166 & 0.0529 \\
    \hline
    \hline
    \multirow{6}[0]{*}{RVR} & \multirow{3}[0]{*}{ conductor} & 0.01\% noise & 0.3992 & 0.4405 & 0.3866 & 0.5726 \\
\cline{3-7}          &       & 0.1\% noise & 0.6590 & 0.5265 & 0.4967 & 0.7483 \\
\cline{3-7}          &       & 1\% noise & 1.1334 & 0.8591 & 0.4805 & 0.7573 \\
\clineB{2-7}{2.5}          & \multirow{3}[0]{*}{ resistor} & 0.01\% noise & 1.2038 & 1.6083 & 0.9931 & 1.0125 \\
\cline{3-7}          &       & 0.1\% noise & 3.1351 & 2.6043 & 1.6627 & 1.4349 \\
\cline{3-7}          &       & 1\% noise & 3.4953 & 2.2805 & 1.4426 & 2.0613 \\
    \hline
    \end{tabular}%
  \label{tab:T6_noisy}%
\end{table}%


\section{Discussion}\label{sec:discussion}
We consider each of the driving questions of the manuscript, in turn, and discuss pros and cons of the reconstruction methods.

\subsection{Q1: How does reconstruction quality from electrode data compare to analytic data?} 
Considering Figure~\ref{fig:T7-compare-rad}, we can see that both analytic and simulated electrode data is reconstructed by Calder\'on's method with a similar smoothing effect.  As one would expect, the smoothing effect would be greater as we increased the mollifying parameter, $t$ and lower, but with more Gibbs phenomenon, as we decrease $t$.

The biggest difference when reconstructing from analytic data versus electrode data is that the Fourier regularization parameter, $T_z$ had to be decreased by half in order to achieve similar reconstructions.  The reason for this can be seen on the right of Figure \ref{fig:T7-compare-rad}, where the Fourier data for the electrode cases deviates greatly from the analytic Fourier data to the right of the vertical dashed line.  This deviation occurs even sooner for the 32 electrode case, which is why we see such an overestimate in the inclusion's conductivity for $L=32$.  Decreasing $T_z$ to 1 would provide a more reasonable reconstruction.

We also compared how each of the reconstruction methods performed in reconstructing target T1 from simulated 128 electrode data, see Figure~\ref{fig:T7-compare}. The evaluation metrics for this are in Table~\ref{tab:T7}. All methods reliably located the center of the inclusion.  We note that from a radially symmetric target, we would expect radially symmetric reconstructions with a perfect center.  We suspect the localization errors here are due to the discretization of a spherical inclusion interpolated to a Cartesian grid on which these metrics were computed. We also note the smoothing effect of CGO methods causing lower RVR values. The CGO methods appear to have better dynamic ranges, more accurately recovering the maximum value of the target.  The LS-based methods overestimate the conductivity but have better MSE and MS-SSIM values.
\subsection{Q2: Can complex-valued admittivity targets be recovered?}
Each of the reconstruction methods was able to recover the main structures in the complex admittivity as shown in Figure~\ref{fig:T3cmplx-compare}.  For the example presented here, both the $\texp$ and Calder\'on reconstructions underestimated conductivity and susceptivity of the heart and, as expected, have smoother reconstructions than both LS-based methods.  
The CGO methods underestimated the dynamic range due to underestimating the conductivity and permittivity of the heart.  The reconstructed permittivity for the lungs is slightly overestimated producing a relative error of $15\%$ and $20\%$ for the heart.  The conductivity of the lungs is approximately $20\%$ relative error with D-bar overestimating and Calder\'on underestimating the conductivity.  Interestingly, the CGO methods  outperformed the LS methods in MS-SSIM.  This could be due to the use of non-optimal parameters in the LS-based methods.  The TV method outperformed the Smooth method in DR, MSE, and MS-SSIM.  Visually, the LS-based methods produced superior susceptivity images for the example considered here.  

\subsection{Q3:  How is reconstruction quality affected by the number of electrodes simulated?}
Each reconstruction method was able to recover the T2-B and T3 targets for each number of electrodes considered ($L=128$, 64, and 32), Figs.~\ref{fig:T3-compare} and \ref{fig:T6-compare}.  As the number of electrodes decreased, the heart and lungs (T2-B) were pushed towards the middle of the domain in both the $\texp$ and Calder\'on reconstructions which is likely due to the reduction in the truncation radius ($T_\xi$ and $T_z$) for the Fourier data.  With $L=32$ electrodes (2\% coverage), the resistive artefact in the center of the $x_1x_3$ plane is due to the joining of the two lungs as shown in the $x_1x_2$ plane.  An investigation of if image quality could be improved for fewer electrodes by increasing the size of the electrodes and thus relative surface coverage (here $L=128$, 64, and 32 electrodes corresponded to 8\%, 4\% and 2\%) is outside the scope of this manuscript.  The resistive and conductive targets in the high-contrast example (Fig.~\ref{fig:T6-compare}) are easily identifiable for all levels of electrodes.

Considering first the heart and lungs T2-B example, with $L=128$ electrodes, the TV method produced the best dynamic range, followed by the CGO methods and then the Smooth LS method.  Reducing to 64 electrodes dropped the dynamic range of Calder\'on's method but it rebounded in the $L=32$ case.  This was a surprise since we would expect poorer reconstructions with fewer electrodes/less data. Further study beyond these two data sets is required to investigate this in more detail.  The dynamic range was poorest for the 32 electrode case with the D-bar method, likely due to the decrease in the truncation radius $T_\xi$.  The LS-based methods produced superior dynamic ranges: TV was best for $128$ and $64$ electrodes, while the Smooth LS was the best for $32$ electrodes.  As to be expected, the LS-based methods best localized the targets and had the highest MS-SSIMs for all electrodes considered.  Note that due to the poor separation of the lungs in the CGO methods for the 32 electrode case, the LE and RVR metrics were not computed.

For the high-contrast example, the reconstructed conductivity values stayed fairly consistent as the number of electrodes decreased from 128 to 32 for each method, respectively.


\subsection{Q4: Can high-contrast targets be recovered?}
We push the limits of the CGO methods by considering an example of fifteen times contrast between conductor and resistor.  All reconstructions clearly show a resistor and conductor for each number of electrodes considered.  However, the LS-based methods more accurately recover the conductivity values of the targets, especially the low-conductivity resistor $\sigma=0.1$~S/m.  While the $\texp$ and Calder\'on methods recover the conductivity of the conductor quite well, they significantly overestimate the conductivity of the resistor, failing to obtain the true dynamic range of the example. 

This is not unexpected, e.g. \cite{Hyvonen2018}.  Linearization based methods assume that the admittivity is a small perturbation from a constant.  This is explicitly included in Calder\'on's method but implicitly included in the $\texp$ approximation 
to the D-bar method.  While the complete D-bar method solves the fully non-linear problem, the $\texp$ approximation considered here is considered a {\it `Born'} approximation as the asymptotic behavior of the CGO solutions is used to compute the non-linear Fourier data $\texp$.  Additionally,  the D-bar method assumes that the non-linear scattering data approaches the linear Fourier data as the magnitude of the non-physical parameter $\zeta$ approaches infinity.  In this work, the magnitude of $\zeta$ was much smaller, following the minimal-zeta approach of \cite{Delbary2014}.  For context, using $T_\xi=16$ corresponded to a maximum value of $\|\zeta\|_2 \approx 11.3$, $T_\xi=12$ had $\|\zeta\|_2 \approx 8.4$, and $T_\xi = 7$ had $\|\zeta\|_2 \approx 4.9$.  Both the D-bar and Calder\'on's method rely on inverting a Fourier transform but in practice are forced to employ a low-pass filtering of the data to deal with the finite nature of solving the problem on a computer as well as for improved stability of the reconstruction as the higher frequencies in the Fourier domain are unstable.  Alternative methods for choosing the non-physical parameter $\zeta$, as well as the full non-linear D-bar method, are outside the scope of this current work.

The LS-based methods overestimated the conductivity of the conductive target but far outperformed the CGO methods in recovering the low conductivity value.  The CGO methods produced the most accurate value for the conductive target for 128 and 64 electrodes; all methods produced quite similar errors in reconstructed high conductivity values at the 32 electrode case.  

\subsection{Q5: What is the effect of noise on the reconstructed admittivity?}
We considered three levels of added relative Gaussian noise: $0.01\%$, $0.1\%$ and $1\%$ in  Figs.~\ref{fig:T3_noisy} and \ref{fig:T6_noisy}.  Calder\'on's method appears most stable producing most easily identifiable images of the targets at $1\%$ noise for the T2-B example. For the T3 example, the Smooth LS method appeared least stable at $1\%$ noise.  The CGO methods handle noise differently than the LS-based methods producing smoother images as the noise level increased, due to the reduction of the truncation radius of the Fourier data, whereas the LS-based methods developed oscillations due to the noise and instability of the EIT problem. 
The reconstructed contrast diminished for the CGO methods, worse for the $\texp$ method, as the level of noise increased, again likely due to the reduced truncation radii $T_\xi$ and $T_z$.  The LS-based methods also suffered diminished accuracy in the dynamic range, especially at the $1\%$ noise level.  All methods overestimated the regional volume ratio (RVR) of the T3 resistor for increased levels of noise with the automated segmentation described above.



\subsection{Comparison of Methods}
Each method appeared to outperform the others depending on the example considered and/or metric used.  Overall, the regularized non-linear LS-based methods routinely outperformed the CGO methods in MS-SSIM, MSE, and LE with the exception of the noisy voltage data where the CGO methods were more robust.  One significant advantage of the CGO methods was {\it computational cost}.  The D-bar method took approximately 5 seconds per data set without leveraging parallelization for the computation of the scattering data (Step 1, \eqref{eq:t-scat-tBIE}) and without utilizing the Fast Fourier Transform in Step~2 \eqref{eq:texp-q} which would give additional speedup.  Calder\'on's method took approximately 50 seconds per data set without parallelization for the computation of $\Fhat$ (Step 1.a., averaging over 30 choices of $a$, \ref{eq:betterFhat}) but approximately 5 seconds following Step 1 with one choice of $a$, \eqref{eq:CaldFhat}, the effect of this on reconstruction quality has not been studied here.  Like the D-bar method, our implementation of Calder\'on's method did not use the Fast Fourier transform in Step~2, \eqref{eq:sigCal_pert}, which would also decrease computation time. In contrast, the computation of reconstructions with the regularized LS-based methods took several hours. 




Although the CGO-based methods did not use the CEM in the solution of the inverse problem, but rather used a gap-model formulation for the ND map, the approach was demonstrated to lead to image quality comparable to the 2D CGO counterparts \cite{Muller2017}, \cite{Hamilton2018_Robust}.  The robustness of the algorithms appears comparable for the algorithms considered based on identifiability of targets across the examples tested.  We note, as did \cite{Delbary2014}, that the $\texp$ reconstruction algorithm is sensitive to the choice of this non-physical parameter $\zeta$.  As $\texp$ depends on both $\zeta$ and $\xi$, the sensitivity in $\zeta$ can be controlled via the low-pass filtering in the $\xi$ grid.  The choice of $T_\xi$ appears to be similar to the 2D case inasmuch as the scattering data visibly `blows up' to $\pm\infty$ outside of a stable ball of radius $T_\xi$. A non-uniform truncation can also be used.  A more thorough study of how best to choose a regularization parameter is left for future work (see \cite{HamiltonLionheartAdler_2019} for suggestions in 2D). Similarly, in Calder\'on's method, the values of $\Fhat(z)$ `blow up' for large $z$ in the presence of noise. In this study, $T_z$ is chosen to minimize that effect, as is the case in 2D. As with the $\texp$ D-bar approach, finding a strategy for choosing the best regularization parameter is the subject of future study. For regularized LS-based methods, there exist approaches for systematically choosing the regularization parameters. For example, the L-curve method \cite{hansen1992analysis}, Morozov's discrepancy principle \cite{blomgren2002modular} or the approach proposed in \cite{gonzalez2016experimental} could be used. However, these approaches are traditionally formulated for a single regularization parameter, and in this manuscript, regularizing the real and imaginary parts of admittivity separately resulted in six regularization parameters for the smoothness regularization and four regularization parameters for the total variation regularization. Due to the high number of parameters to tune, in this manuscript, none of the systematic approaches were used and the regularization parameters were chosen manually based on visual inspection of the reconstructions.
\section{Conclusions}\label{sec:conclusion}
This manuscript presented the first reconstructions of CGO methods on a sphere with CEM electrode data, and provides the first direct comparison of CGO methods with common regularized non-linear least-squares methods in three-dimensions.  A variety of targets were considered, including high-contrast and complex admittivities. The effect of noisy voltage data, as well varying the number of electrodes, were explored.  More complicated domain shapes and experimental tank studies are left for future work.

We have demonstrated that the CGO methods are now a viable option for realistic data.  The CGO methods appear to have comparable robustness to noise and may outperform the LS-based methods for higher levels of noise (targets are identifiable even if dynamic range is lower).  The speed and reliability of the CGO methods to routinely produce useful information about the targets indicates they could be a viable method for initializing common optimization-based algorithms.   Due to the high computational cost of 3D optimization-based methods, a better initial guess could provide a significant reduction in the number of iterations required.  For applications that only require the detection of a conductor or resistor, these CGO methods provide a standalone fast means of identification.


\section*{Acknowledgments}
SH and DI were supported by the National Institute Of Biomedical Imaging And Bioengineering of the National Institutes of Health under Award Number R21EB028064. The content is solely the responsibility of the authors and does not necessarily represent the official views of the National Institutes of Health.  
JT and VK were supported by the Academy of Finland (Project 312343, Finnish Centre of Excellence in Inverse Modelling and Imaging), the Jane and Aatos Erkko Foundation and Neurocenter Finland.  


\bibliographystyle{amsalpha}
\small
\bibliography{bibliographyRefs}
\end{document}